\documentclass[english,aps,pre,twocolumn,showpacs,floatfix]{revtex4-1}
\usepackage{amsmath}
\usepackage{graphicx}
\usepackage{amssymb}
\usepackage{dcolumn} 
\usepackage{bm} 
\usepackage{color}
\usepackage{subfigure}
\usepackage[none]{hyphenat}

\makeatletter

\usepackage{amsfonts}
\usepackage{epsfig}
\usepackage{dsfont}
\usepackage{babel}

\makeatother

\begin{document}
\title{Method for generating two coupled Gaussian stochastic processes}

\author{Tayeb \surname {Jamali}}
\email[Email: ]{jamalitayeb@gmail.com}
\affiliation{Department of Physics, Shahid Beheshti University, G.C., Evin, Tehran 19839, Iran}

\author{G. R. \surname {Jafari}}
\affiliation{Department of Physics, Shahid Beheshti University, G.C., Evin, Tehran 19839, Iran}

\date{\today}

\begin{abstract}
Most processes in nature are coupled; however, extensive null models for generating such processes still lacks. We present a new method to generate two \emph{coupled} Gaussian stochastic processes with arbitrary correlation functions. This method is developed by modifying the Fourier filtering method. The robustness of this method is proved by generating two coupled fractional Brownian motions and extending its range of application to Gaussian random fields.

\end{abstract}


\maketitle

\section{Introduction}
\label{sec: intro}

Evolution of real systems mainly contains correlated noises in which some valuable information may be hidden. To study such noises numerically, a number of methods have been developed such as the matrix decomposition methods~\cite{HJ2013,GL1996}, the Fourier filtering method (FFM)~\cite{SFI,PHSS91,PHSS92,MHSS96}, and the circulant embedding method~\cite{WC1994,DN97,LPS2014}.  All these methods are capable of generating realizations of only one stochastic process. Therefore, when facing multiple processes simultaneously, these methods can only be used to the individual study of each process. But, since the evolution of a process in reality depends more or less on the evolution of other processes, an individual study of a process would lead to vague results.

In order to measure the coupling between two or more processes, some various techniques have been introduced, e.g. the detrended cross-correlation analysis (DXA)~\cite{ps,Shadkhoo,DXA2014}, partial-DXA~\cite{partialDXA}, coupled-DXA~\cite{HVJ}, the cross-wavelet analysis~\cite{wavelet,djkms}, the random matrix theory~\cite{lcbp,pgras2,pgrags}, the coupled level-crossing method~\cite{JJF}, the cross-visibility algorithm~\cite{Mehraban} and etc. There also exist special methods for generating two coupled processes but only for power-law correlations~\cite{PFSI,PHLEI}. However, there is still lack of a general method.

In this paper, we aim to present a method capable of generating two Gaussian coupled stochastic processes with any arbitrary correlations. To do this, we modify the Fourier filtering method. The established FFM can be thought of as a machine which takes an uncorrelated sequence of random numbers and returns a correlated sequence of random numbers. We modify this machinery so that it takes two uncorrelated sequences of random numbers and returns two correlated ones. To be more precise, linear combinations of two Gaussian white noises in the Fourier space is written, such that when returning to the real space, our desired correlations would be injected in the processes. Then, the method is implemented for modeling two coupled Brownian motions, and two coupled fractional Brownian motions (FBM). We also show how to extend this method to generate two coupled Gaussian random fields.

\section{Method}
\label{sec:method}
Consider two independent sequences of $L$ uncorrelated random numbers $\{u_i\}_{i=1}^L$ and $\{v_i\}_{i=1}^L$, with Gaussian distributions. Their correlation functions are $\langle u_i u_{i+n} \rangle\sim\delta_{n,0}$ and $\langle v_i v_{i+n} \rangle\sim\delta_{n,0}$, where $\delta_{n,0}$ is the Kronecker delta, and $\langle\dots\rangle$ denotes an ensemble average. The aim of this study is to provide an algorithm by which one could construct two stationary sequences $\{x_i\}$ and $\{y_i\}$  with desired autocorrelations
\begin{subequations}
\begin{align}
\label{eq: ACFx}
C_{xx}(n) &=\langle x_i x_{i+n} \rangle \\
\label{eq: ACFy}
C_{yy}(n) &=\langle y_i y_{i+n} \rangle,
\end{align}
and the cross-correlation
\begin{align}
C_{xy}(n) &=\langle x_i y_{i+n} \rangle, \label{eq: CCF}
\end{align}
\end{subequations}
{starting from the two sequences $\{u_i\}$ and $\{v_i\}$. In this line, we take the advantage of Fourier filtering method for generating two coupled processes, $\{x_i\}$ and $\{y_i\}$. Working in the Fourier space brings the need to deal with spectral densities instead of correlation functions. Here, the corresponding autospectral densities are
\begin{subequations}
\begin{align}
\label{eq: SDxx}
S_{xx}(q) &=\langle x_q x_{-q} \rangle \\
\label{eq: SDyy}
S_{yy}(q) &=\langle y_q y_{-q} \rangle,
\end{align}
and a cross-spectral density is
\begin{align}
\label{eq: SDxy}
S_{xy}(q) &=\langle x_q y_{-q} \rangle,
\end{align}
\end{subequations}
where $\{x_q\}$ and $\{y_q\}$ are Fourier transforms of $\{x_i\}$ and $\{y_i\}$ respectively. The autospectral densities $S_{xx}(q)$ and $S_{yy}(q)$ are real-valued even functions while the cross-spectral density $S_{xy}(q)$ is a complex-valued function~\cite{Bendat}. It is noteworthy to recall that for a stationary random process, due to the Wiener-Khintchine theorem, the spectral density is the Fourier transform of the correlation function~\cite{Wiener}.

In order to construct two stationary sequences $\{x_i\}$ and $\{y_i\}$ with the desired correlations we suggest combining two independent and uncorrelated sequences $\{u_i\}$ and $\{v_i\}$ in the Fourier space as
\begin{equation}
\label{eq: xq and yq}
\begin{aligned}
x_q &=A_q u_q+B_q v_q \\
y_q &=C_q u_q+D_q v_q,
\end{aligned}
\end{equation}\\
where $\{u_q\}$ and $\{v_q\}$ are the Fourier transforms of $\{u_i\}$ and $\{v_i\}$ respectively. The coefficients $A_q, B_q, C_q$, and $D_q$ are functions in the Fourier space which are of essential importance for specifying $x_q$ and $y_q$. The coefficients $A_q, B_q, C_q$, and $D_q$ are well named as being of essential importance because having these coefficients would lead us to obtain two sequences $\{x_i\}$ and $\{y_i\}$ possessing the desired correlations. In this regard, after substituting $x_q$ and $y_q$ from Eq.~(\ref{eq: xq and yq}) into the expressions for spectral densities (Eqs.~(\ref{eq: SDxx})-(\ref{eq: SDxy})), the relation of these coefficients with the desired spectral densities is obtained
\begin{subequations}
\begin{align}
\label{eq: SDxx and Coefficients}
S_{xx}(q) &=A_q A_{-q}+B_q B_{-q} \\
\label{eq: SDyy and Coefficients}
S_{yy}(q) &=C_q C_{-q}+D_q D_{-q} \\
\label{eq: SDxy and Coefficients}
S_{xy}(q) &=A_q C_{-q}+B_q D_{-q}.
\end{align}
\end{subequations}
To derive these relations we use the fact that the spectral density of a white noise is constant~\cite{Bendat}. From all possible answers for the coefficients we only take into account the case
\begin{equation}
\label{eq: ABCD}
\begin{aligned}
A_{-q} &=A_q^* ,\qquad B_{-q}=B_q^* \\
C_{-q} &=C_q^* ,\qquad D_{-q}=D_q^*.
\end{aligned}
\end{equation}
This exactly resembles the characteristics of the spectral densities~\cite{Bendat}. This assumption does not result in any loss of information, because in the end, we find a class of coefficients by which two stationary sequences with desired correlations can be generated. The Eq.~(\ref{eq: ABCD}) leads to
\begin{subequations}
 \begin{align}
 \label{eq: SDxx and Coefficients 2}
   S_{xx}(q) &= |A_q|^2 + |B_q|^2 \\
\label{eq: SDyy and Coefficients 2}
   S_{yy}(q) &= |C_q|^2 + |D_q|^2 \\
\label{eq: SDxy and Coefficients 2}
   S_{xy}(q) &= A_q C_q^* + B_q D_q^*.
  \end{align}
\end{subequations}
Before solving this system of equations it is instructive to interpret it algebraically. Equations.~(\ref{eq: SDxx and Coefficients 2}) and~(\ref{eq: SDyy and Coefficients 2})  are nothing but the square of Euclidean length of two complex vectors $\mathbf{\Phi}\equiv(A_q,B_q)$ and $\mathbf{\Psi}\equiv(C_q,D_q)$. Equation.~(\ref{eq: SDxy and Coefficients 2}) also represents the scalar product $\langle \mathbf{\Psi},\mathbf{\Phi}\rangle$. Hence, Eqs.~(\ref{eq: SDxx and Coefficients 2})-(\ref{eq: SDxy and Coefficients 2}) talk about the length of two complex vectors and the angle between them~\cite{LF,Scharnhorst}. This algebraic interpretation leads us to the following answer to the system of equations ~(\ref{eq: SDxx and Coefficients 2})-(\ref{eq: SDxy and Coefficients 2}) which also satisfies the condition of Eq.~(\ref{eq: ABCD})
\begin{equation}
\label{eq: Coefficients 1}
\begin{aligned}
A_q &= \sqrt{S_{xx}(q)} \cos \alpha_q \quad,\quad B_q = \sqrt{S_{xx}(q)} \sin \alpha_q \\
C_q &= \sqrt{S_{yy}(q)} \cos \beta_q \quad,\quad D_q = \sqrt{S_{yy}(q)} \sin \beta_q,
\end{aligned}
\end{equation}
where $\alpha_q$ and $\beta_q$ should satisfy
\begin{equation}
\label{eq: Coefficients 2}
\alpha_q-\beta_q^* = \arccos\left(\frac{S_{xy}(q)}{\sqrt{S_{xx}(q)S_{yy}(q)}}\right),
\end{equation}
and $\beta_q^*$ is the complex conjugate of $\beta_q$.

\section{Algorithm}
\label{sec: algorithm}
We propose a numerical algorithm using the Eqs.~(\ref{eq: Coefficients 1}) and (\ref{eq: Coefficients 2}) in order to generate two stationary sequences with the desired correlations. The algorithm is as follows:\\
a) Generate two independent sequences of uncorrelated random numbers  $\{u_i\}$ and $\{v_i\}$ with a Gaussian distribution, then calculate their Fourier transform coefficients $\{u_q\}$ and $\{v_q\}$.\\
b) Calculate the fourier transforms $S_{xx}$, $S_{yy}$ and $S_{xy}$ of the desired correlation functions $C_{xx}$, $C_{yy}$ and $C_{xy}$. \\
c) Obtain the coefficients $A_q, B_q, C_q$ and $D_q$ using Eq.~(\ref{eq: Coefficients 1}) \& Eq.~(\ref{eq: Coefficients 2}) and substitute in Eq.~(\ref{eq: xq and yq}) to get $\{x_q\}$ and $\{y_q\}$. \\
d) Calculate the inverse Fourier transform of $\{x_q\}$ and $\{y_q\}$ in order to obtain two sequences  $\{x_i\}$ and $\{y_i\}$ with the desired correlations.

\section{Applications}
\label{sec: applications}

\subsection{Two coupled Brownian motions}
\label{sec: Two coupled Brownian motions}

The method proposed in the present study provides two Gaussian time series with any sort of autocorrelation and cross-correlation. For instance, this method enables the construction of two coupled time series without any autocorrelations. In other words, it allows us to produce two Gaussian white noises $\{x_i\}$ and $\{y_i\}$ which are coupled. If we consider $\{x_i\}$ and $\{y_i\}$ as the steps of two random walks, then, $X(t)=\sum_{i=1}^t x_i$ and $Y(t)=\sum_{i=1}^t y_i$ would represent the positions of two coupled Brownian motions at time $t$. In Fig.~\ref{fig: coupled_Bm}, two coupled Brownian motions are shown for three different kind of couplings. The top panel of Fig.~\ref{fig: coupled_Bm} refers to a Gaussian coupling, the middle panel corresponds to an exponential coupling, and the bottom panel shows a damped harmonic coupling. Note that in generating the two series, $\{x_i\}$ and $\{y_i\}$, the sequences $\{u_i\}$ and $\{v_i\}$ have been considered initially equal for all three panels because the coupling shows its effects in the pattern of two series $X(t)$ and $Y(t)$ in the right panels.

\begin{figure}[tb]
\centering
\includegraphics[width=1.0\linewidth]{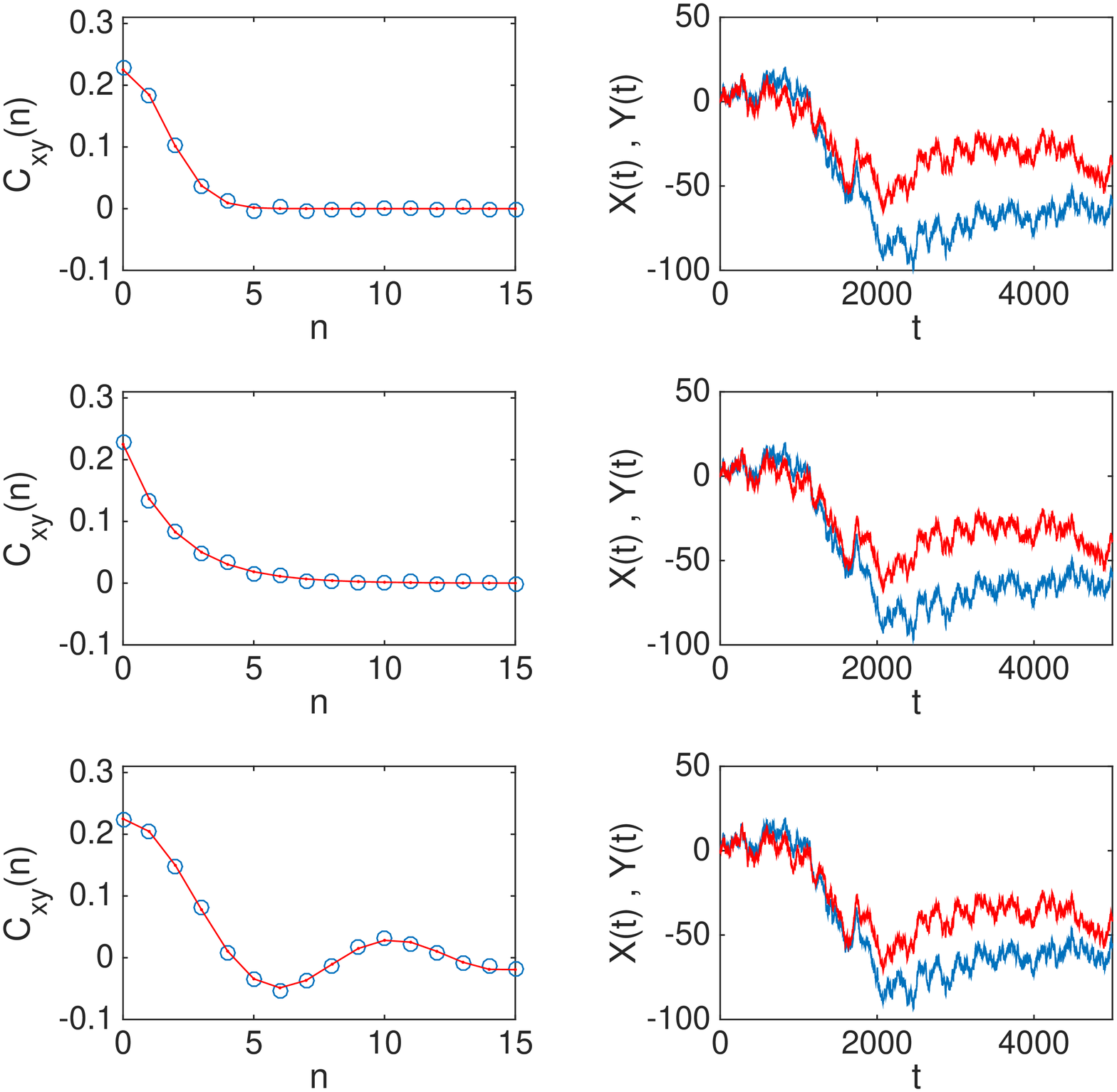}
\caption{Left panels show the coupling between two white noises $\{x_i\}$ and $\{y_i\}$ of the length $L=2^{10}$ and the right panels show the realizations of two coupled Brownian motions $X(t)=\sum_{i=1}^t x_i$ and $Y(t)=\sum_{i=1}^t y_i$. Top panels are related to Gaussian coupling, middle panels is for exponential coupling, and bottom panels represent damped harmonic coupling. Note that the circles show couplings obtained numerically from our algorithm which very well coincided with the expected couplings shown by the solid lines.}
\label{fig: coupled_Bm}
\end{figure}

\begin{figure}[t!]
\centering
    \subfigure{
        \centering
        \includegraphics[width=0.23\textwidth]{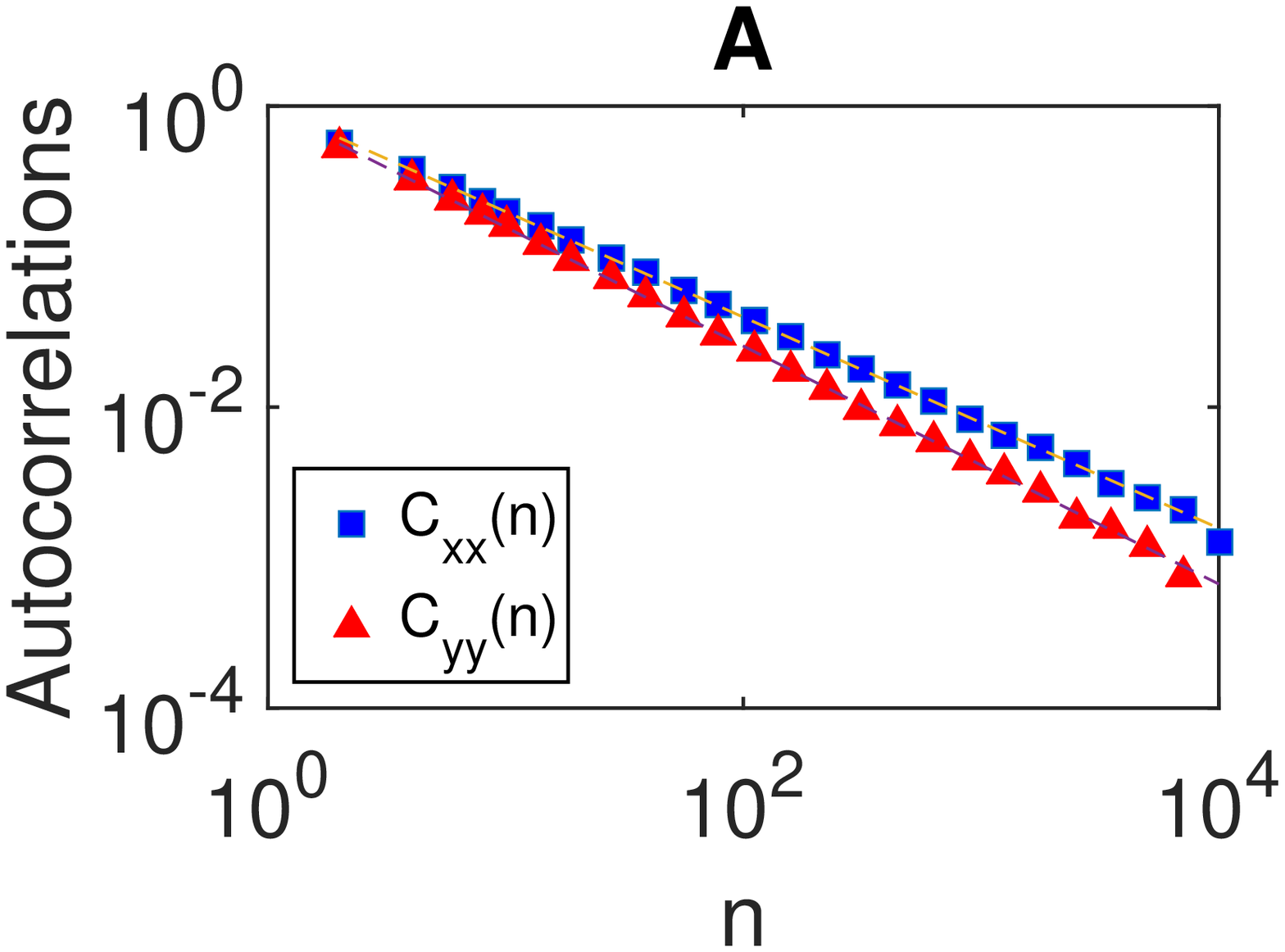}
        \label{fig: blue community}
    }
    \hspace*{-0.9em}
    \subfigure{
        \centering
        \includegraphics[width=0.23\textwidth]{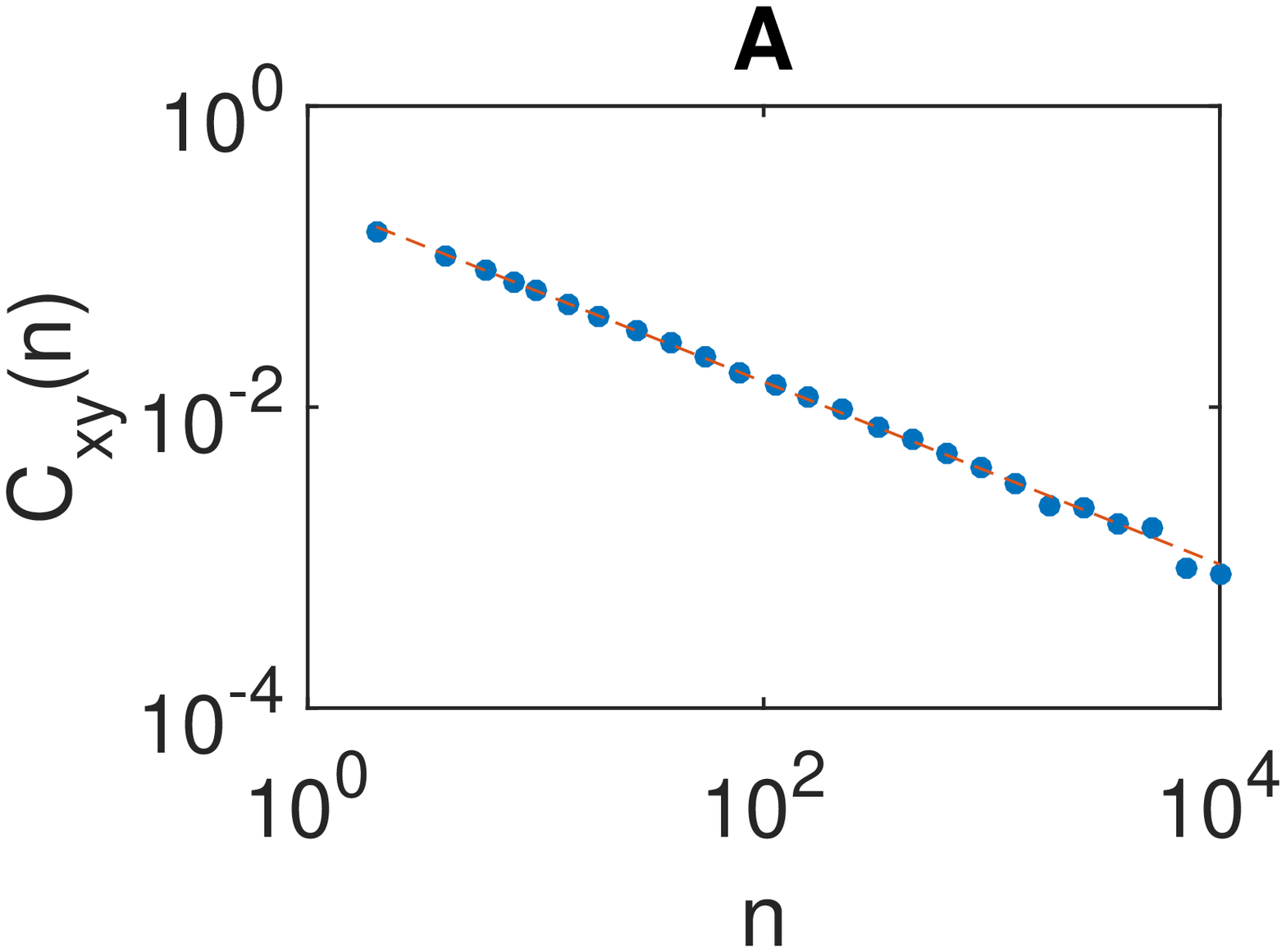}
        \label{fig:subfigure2}
    }

    \subfigure{
        \centering
        \includegraphics[width=0.23\textwidth]{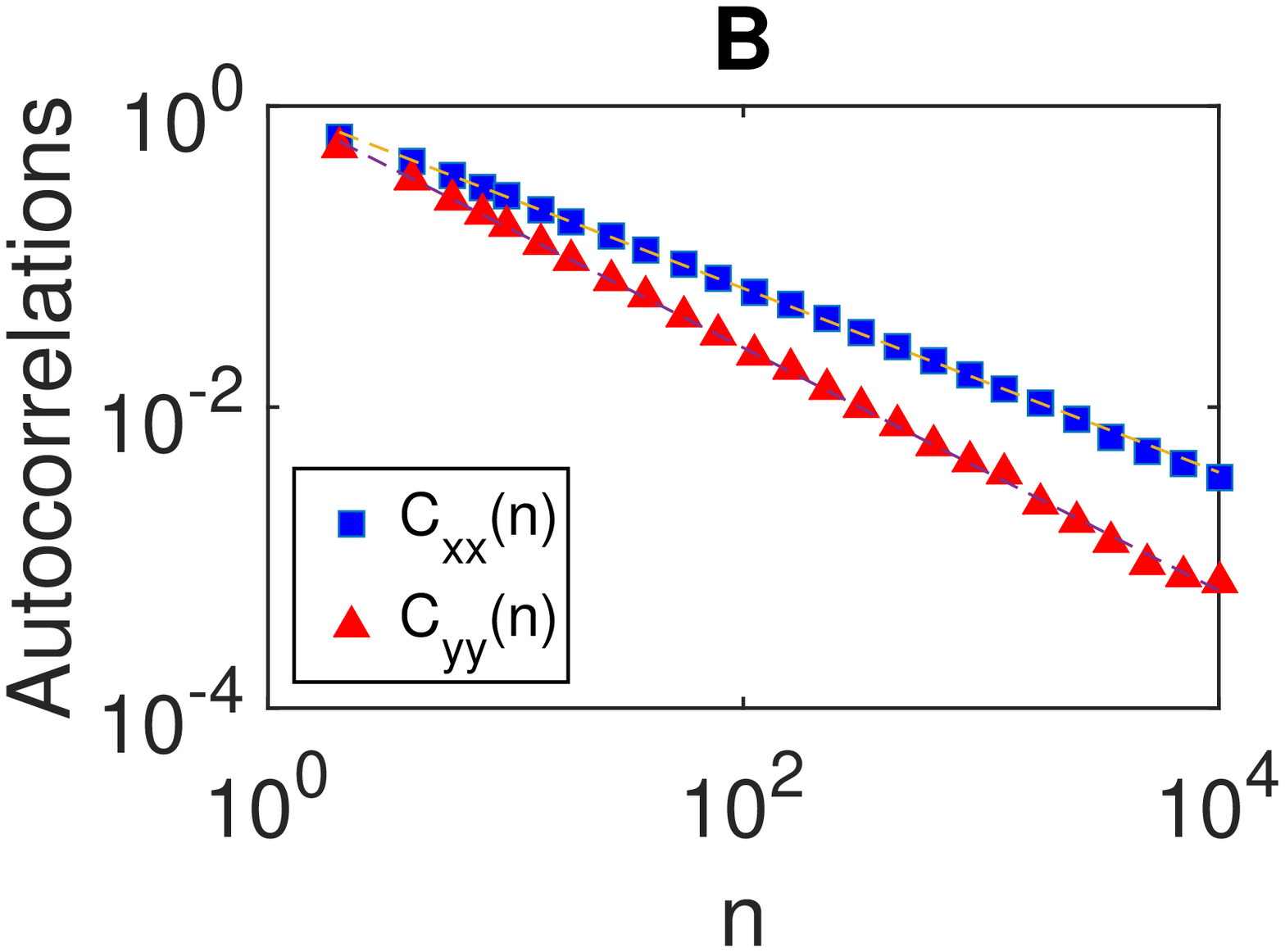}
        \label{fig:subfigure2}
    }
    \hspace*{-0.9em}
    \subfigure{
        \centering
        \includegraphics[width=0.23\textwidth]{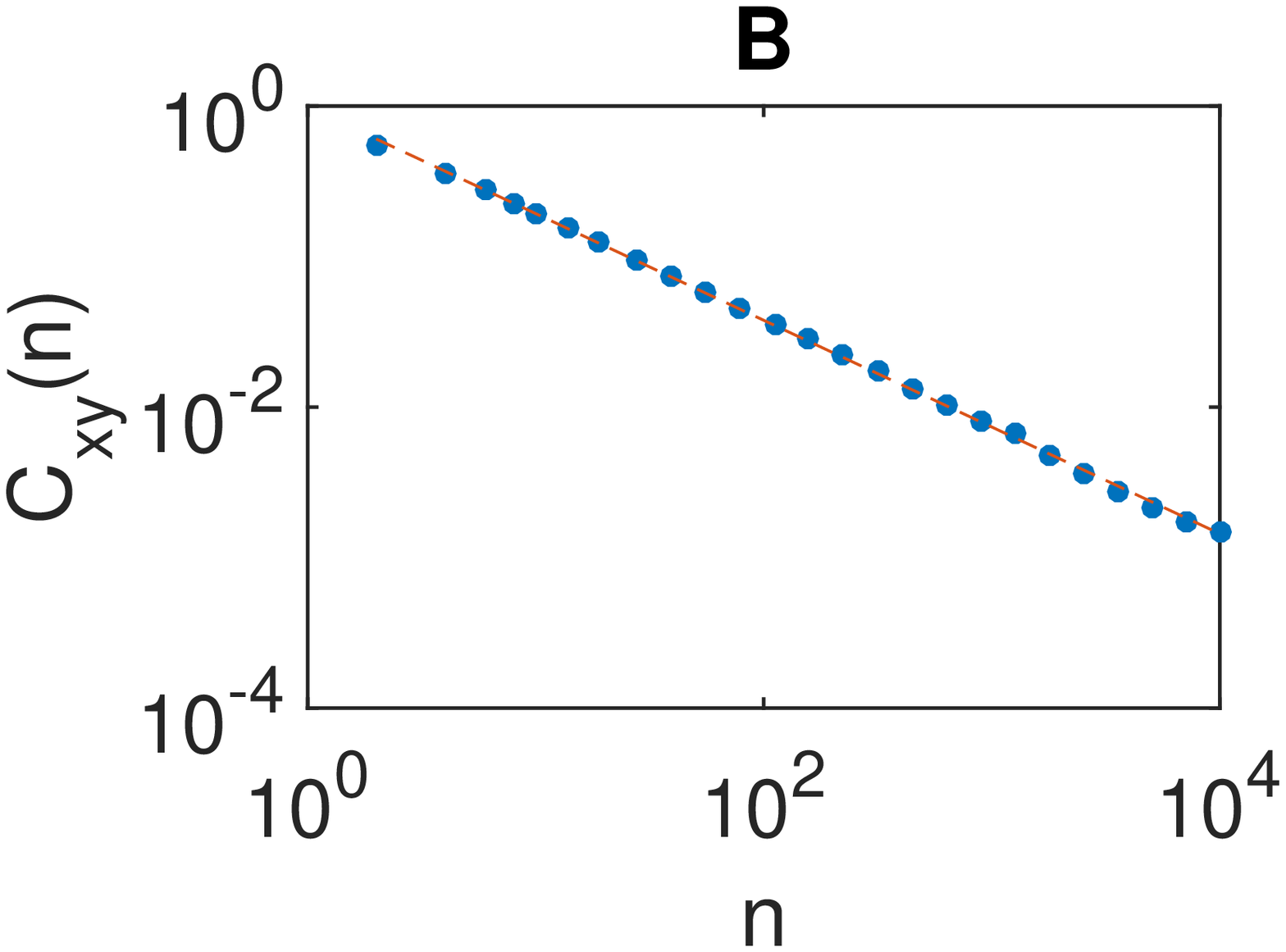}
        \label{fig:subfigure3}
    }

    \subfigure{
        \centering
        \includegraphics[width=0.23\textwidth]{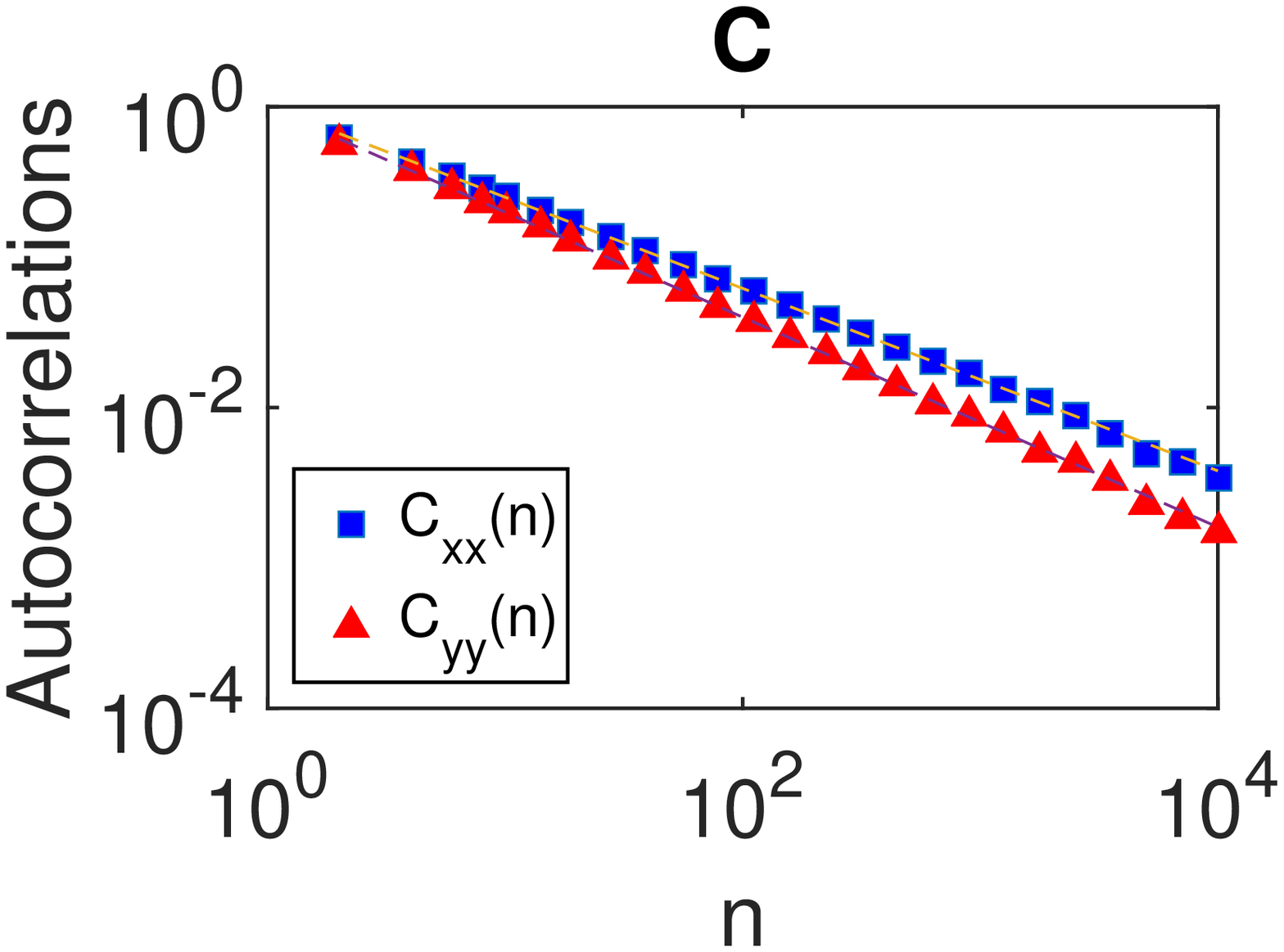}
        \label{fig: red community}
    }
    \hspace*{-0.9em}
    \subfigure{
        \centering
        \includegraphics[width=0.23\textwidth]{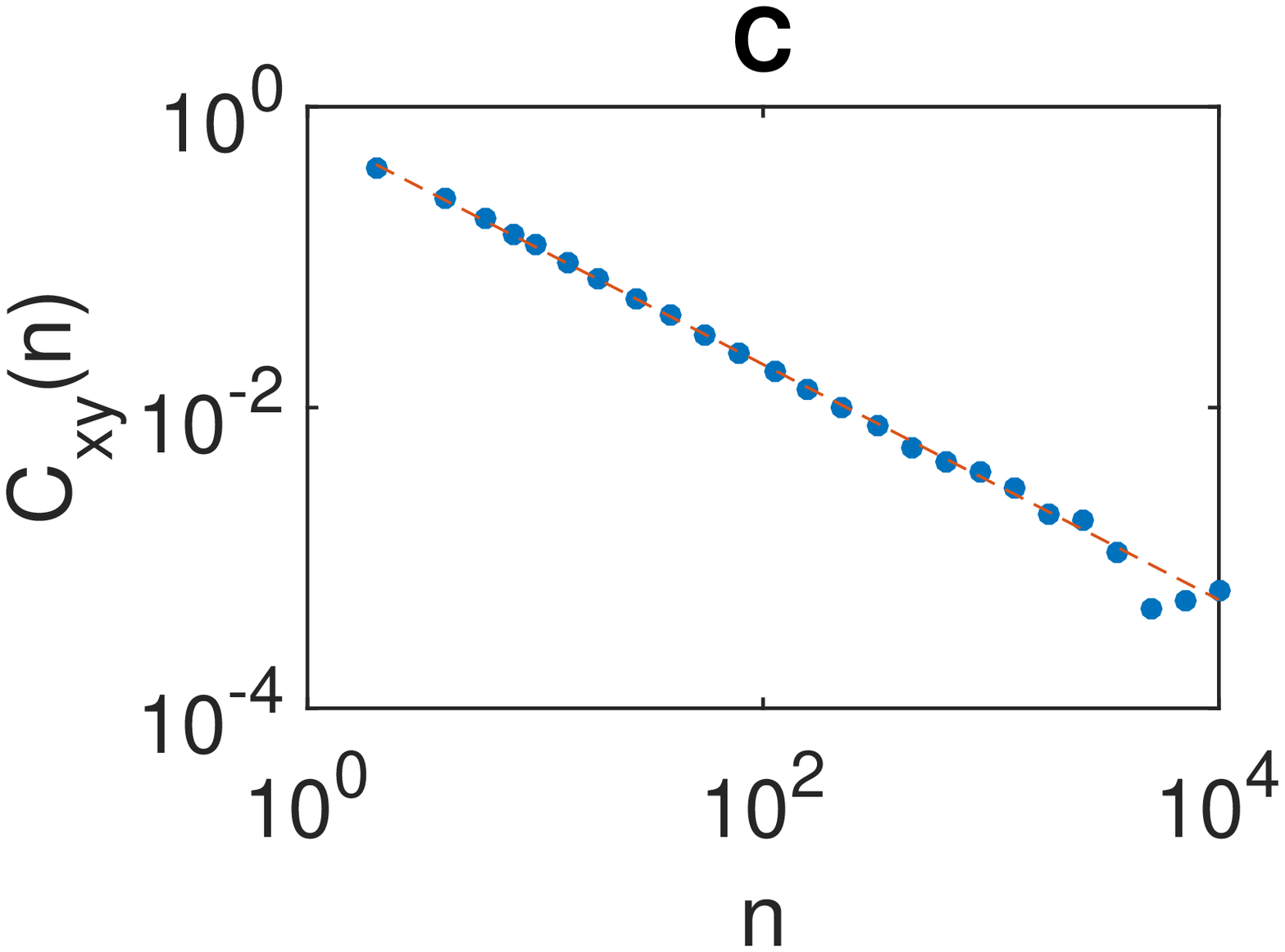}
        \label{fig: red community}
    }
\caption{The log-log plots of the average correlations $C_{xx}(n)$, $C_{yy}(n)$, and $C_{xy}(n)$ for three different cases labeled by 'A', 'B', and 'C'. Each case is the average of $100$ power-law correlated samples for $L=2^{21}$. For cases 'A', 'B', and 'C' the desired exponents $(\gamma_{xx},\gamma_{yy},\gamma_{xy})$ are respectively equal to $(0.7,0.8,0.6)$, $(0.6,0.8,0.7)$, and $(0.6,0.7,0.8)$. The dashed lines represent the best fits for the exponents $(\gamma_{xx},\gamma_{yy},\gamma_{xy})$ which yeald values of $(0.70\pm 0.01,0.80\pm 0.01,0.60\pm 0.03)$, $(0.61\pm 0.01,0.80\pm 0.02,0.71\pm 0.02)$ and $(0.60\pm 0.01,0.70\pm0.01,0.79\pm 0.05)$ for cases 'A', 'B', and 'C' respectively.}
\label{fig: coupled_fBm}
\end{figure}

\subsection{Two coupled fractional Brownian motions}
\label{sec: Two coupled fractional Brownian motions}

Although there are many techniques developed to generate a time series with a power-law autocorrelation~\cite{MHSS95,MHSS96,HS,VJ}, but in reality we sometimes face a couple of power-law autocorrelated time series which posses a power-law cross-correlation e.g. fluctuations of stock prices~\cite{PGGAS}. Podobnik et. al. modeled two power-law autocorrelated time series with long-range cross-correlations~\cite{PFSI,PHLEI}. In this section we aim to generate such power-law correlated time series by our algorithm where we use the special form of correlation function proposed by Makse et. al.~\cite{MHSS96}
\begin{equation}
\label{eq: power-law correlation function}
C(l)=(1+l^2)^{-\gamma/2},
\end{equation}
where $0<\gamma<1$ is the correlation exponent. The function $C(l)$ which is well defined at $l=0$ shows the desired power-law behavior for large $l$. This special form for power-law correlation function has two advantages. First, it enables generating a power-law correlation throughout the system. Second, its Fourier transform, spectral density function,  has an analytic form stated by
\begin{equation}
\label{eq: power-law spectral density}
S(q)=\frac{2\pi^{1/2}}{\Gamma(\beta+1)}\left(\frac{q}{2}\right)^\beta K_\beta(q),
\end{equation}
where $K_\beta(q)$ is the modified Bessel function of order $\beta=(\gamma-1)/2$, and $\Gamma$ is the gamma function, for more details see ref~\cite{MHSS96}. Here, we use Eq.~(\ref{eq: power-law spectral density}) for spectral densities $S_{xx}, S_{yy},$ and $S_{xy}$ with exponents $\gamma_{xx}, \gamma_{yy},$ and $\gamma_{xy}$ respectively. As an example, we applied our algorithm to generate two power-law autocorrelated series with a power-law cross-correlation for three different classes of exponents $\gamma_{xx}, \gamma_{yy},$ and $\gamma_{xy}$. In each row of Fig.~\ref{fig: coupled_fBm} the autocorrelations $C_{xx}$ \& $C_{yy}$ and a cross-correlation $C_{xy}$ is illustrated for 100 samples of the length $2^{21}$.  Note that the values for $\gamma_{xx}, \gamma_{yy},$ and $\gamma_{xy}$ are obtained according to the best fit for every curve in the panels of Fig.~\ref{fig: coupled_fBm}.

Similarly to section~\ref{sec: Two coupled Brownian motions}, we can consider the series $\{x_i\}$ and $\{y_i\}$ which are power law correlated as the steps of two random walks. $X(t)=\sum_{i=1}^t x_i$ and $Y(t)=\sum_{i=1}^t y_i$ represent the positions of two fractional Brownian motions at time $t$ with a power-law coupling.

\begin{figure*}[t!]
\centering
    \subfigure{
        \centering
        \includegraphics[width=0.48\textwidth]{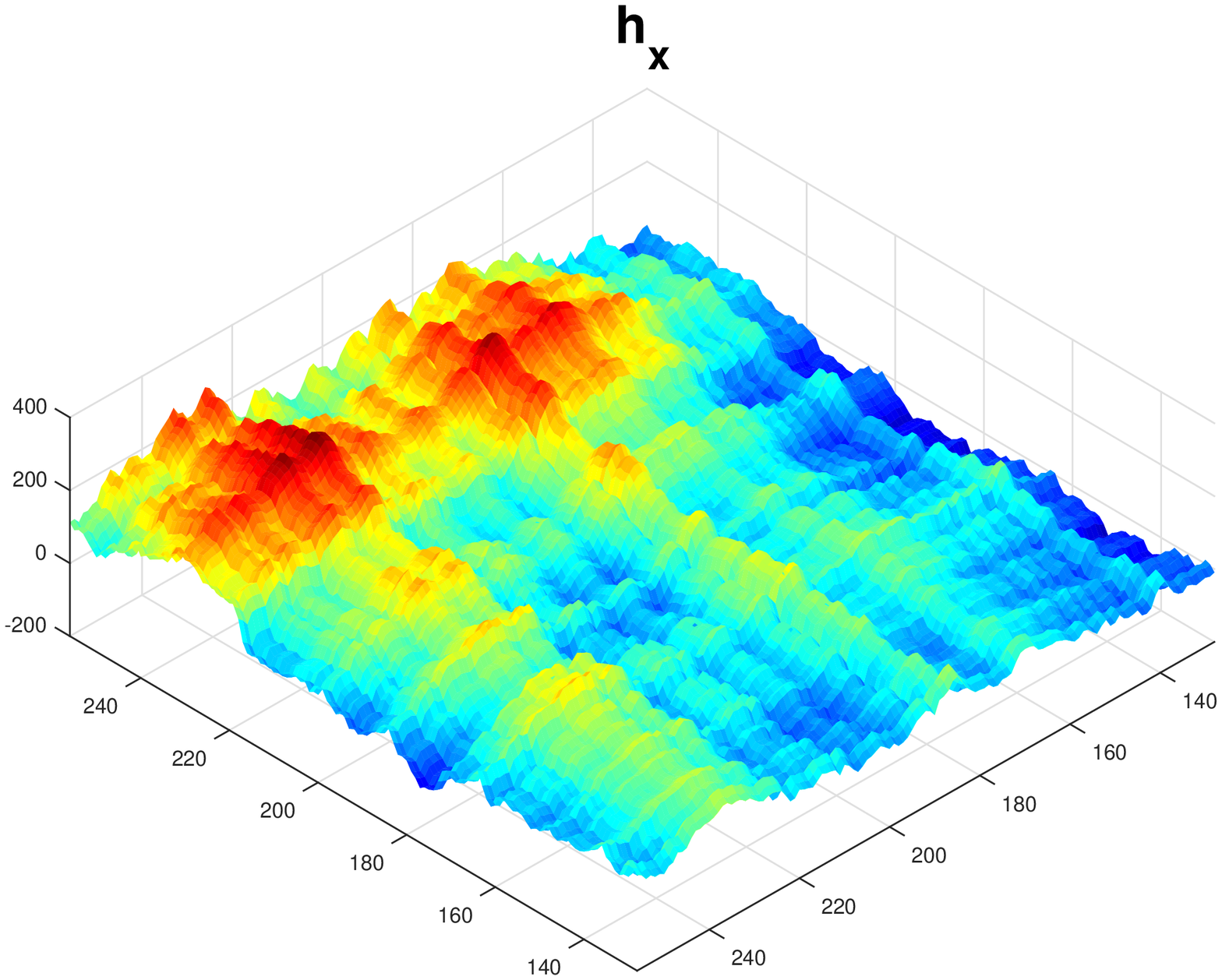}
    }
    \hspace*{-0.9em}
    \subfigure{
        \centering
        \includegraphics[width=0.48\textwidth]{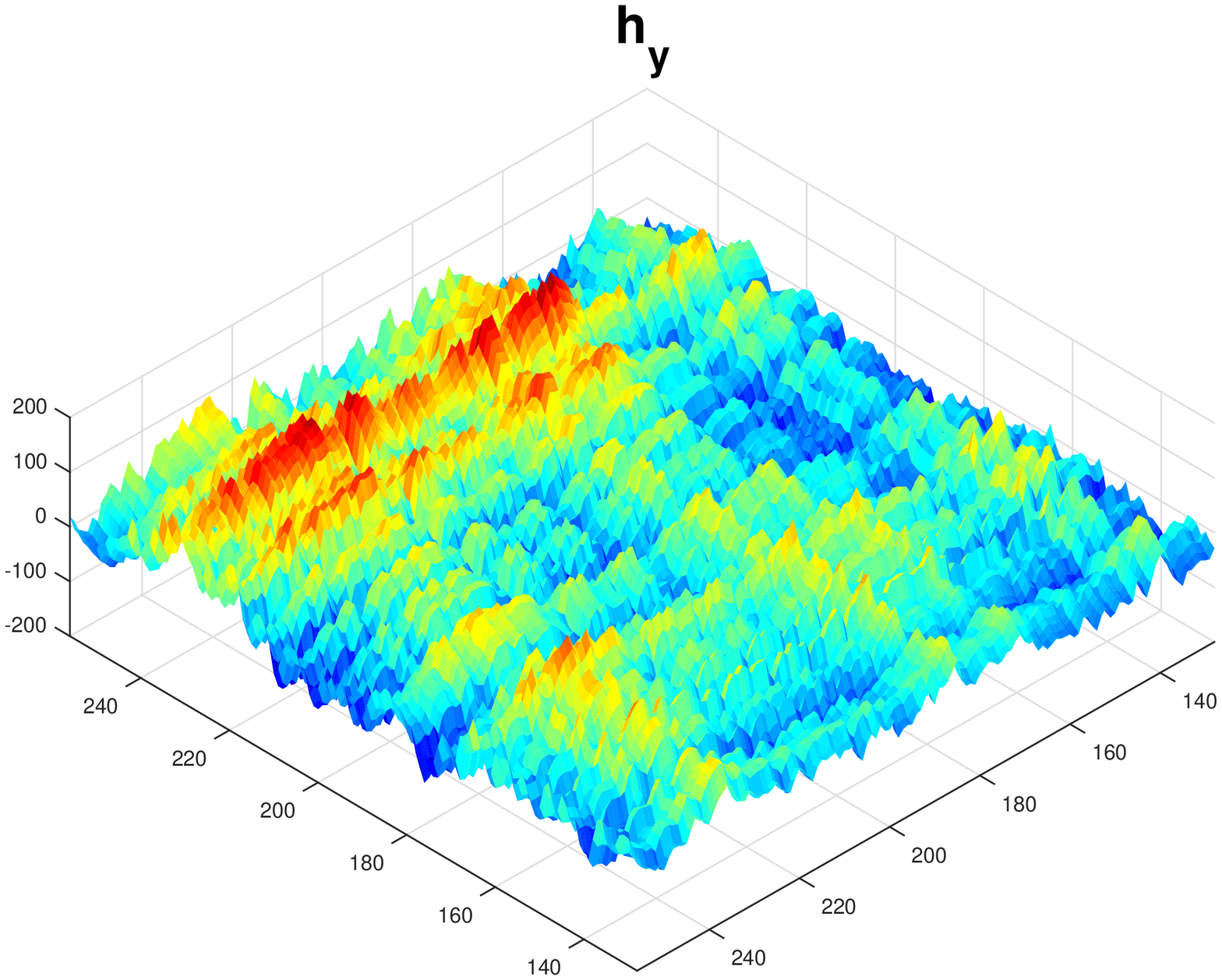}
    }

    \subfigure{
        \centering
        \includegraphics[width=0.32\textwidth]{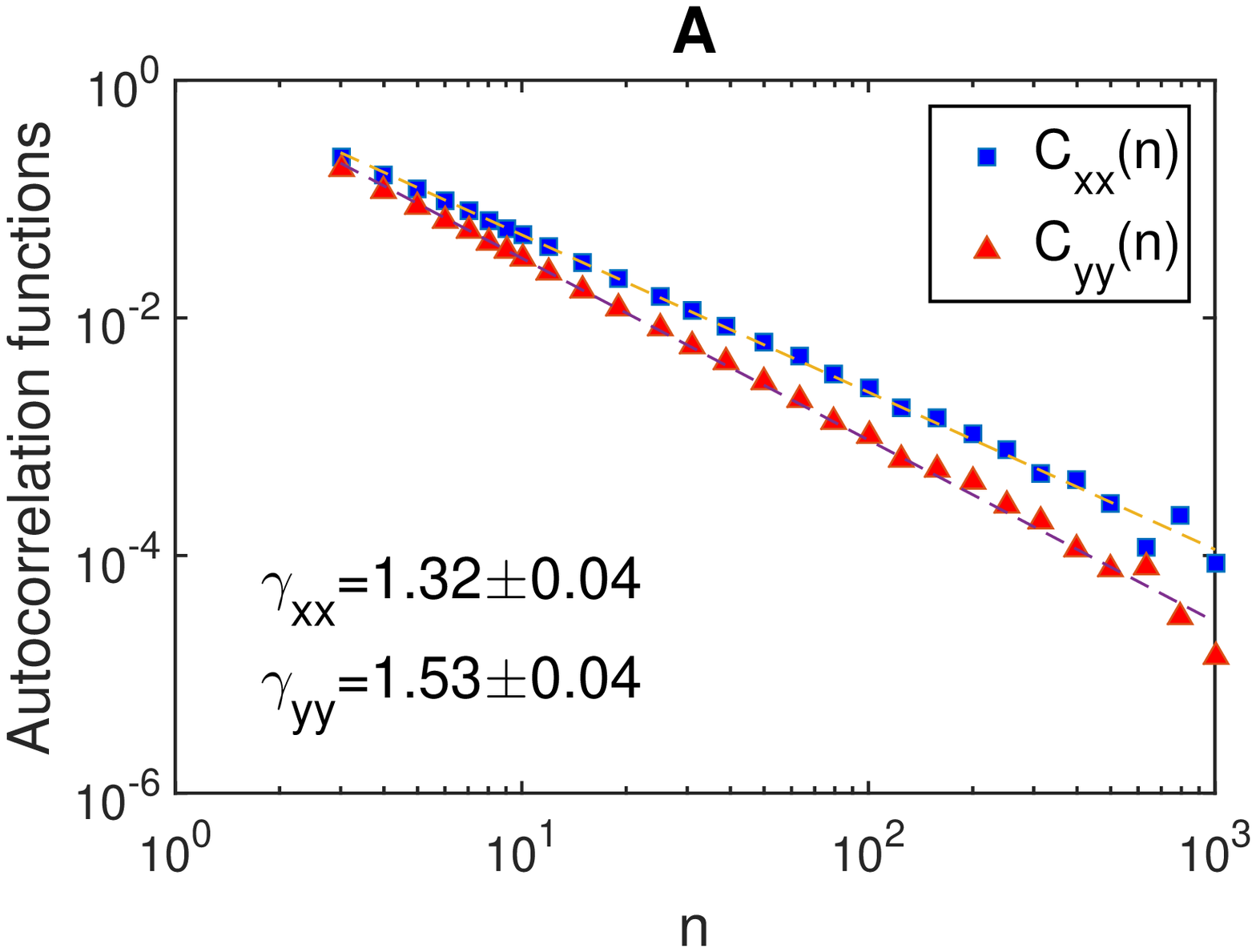}
    }
    \hspace*{-0.9em}
    \subfigure{
        \centering
        \includegraphics[width=0.32\textwidth]{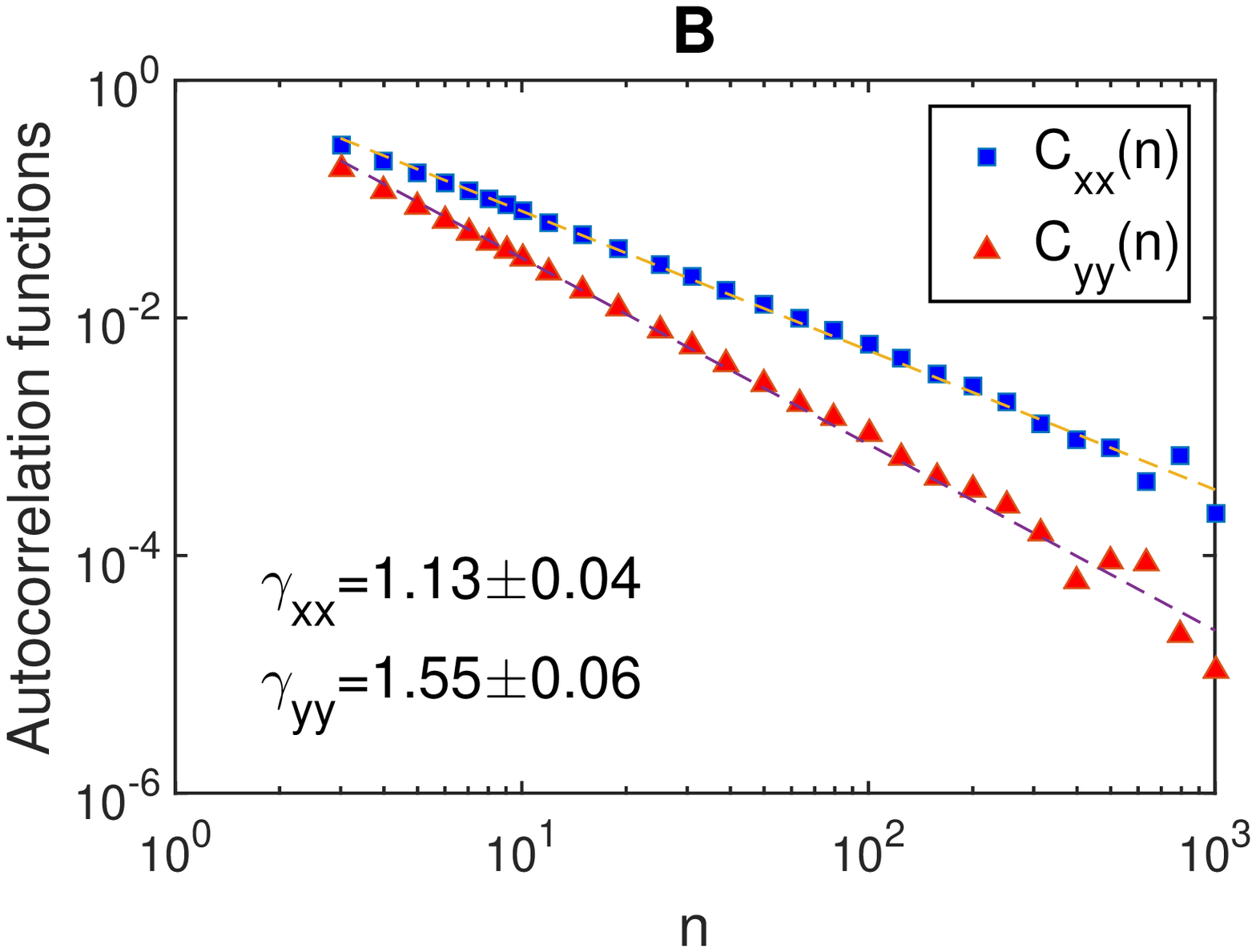}
    }
    \hspace*{-0.9em}
    \subfigure{
        \centering
        \includegraphics[width=0.32\textwidth]{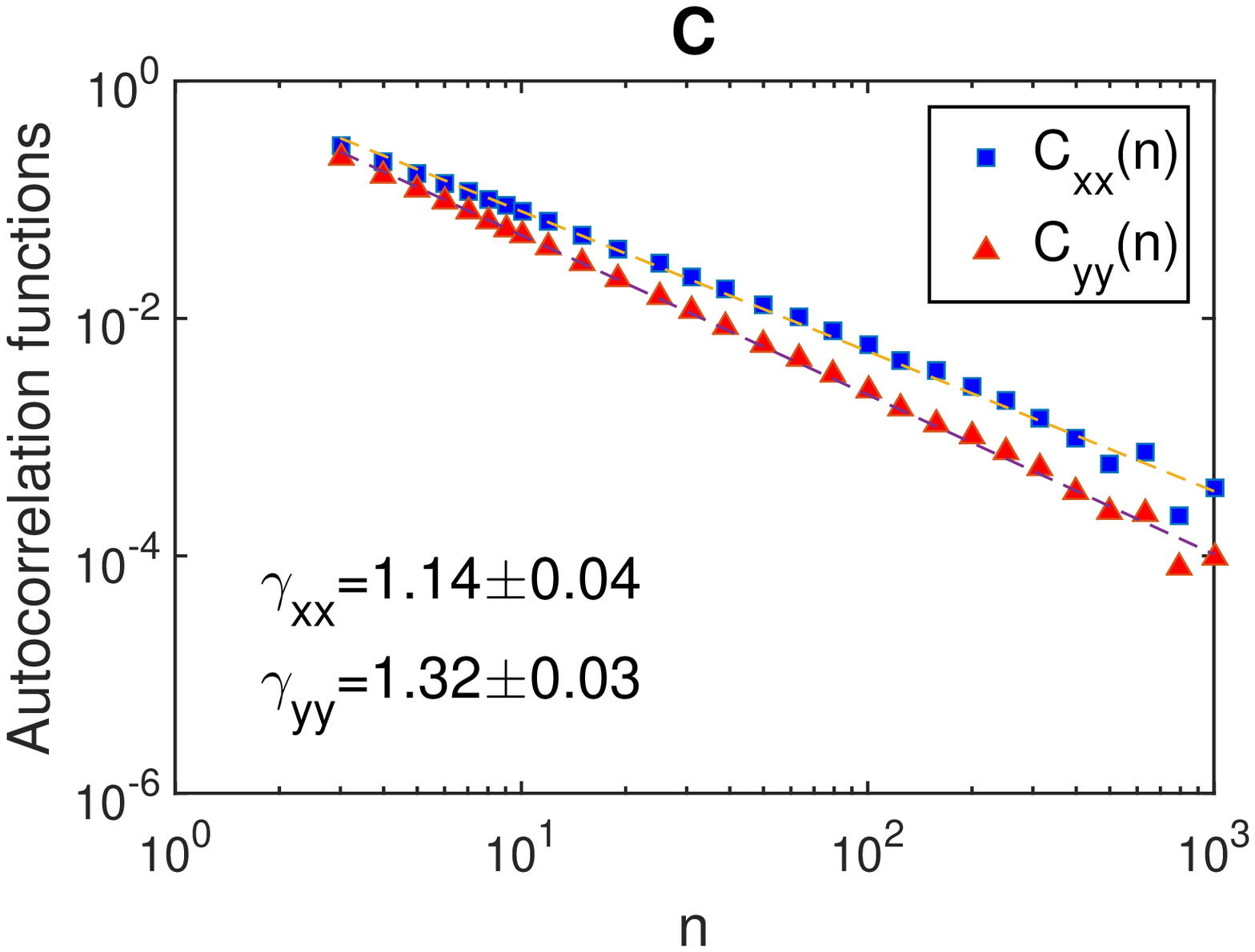}
    }

    \subfigure{
        \centering
        \includegraphics[width=0.32\textwidth]{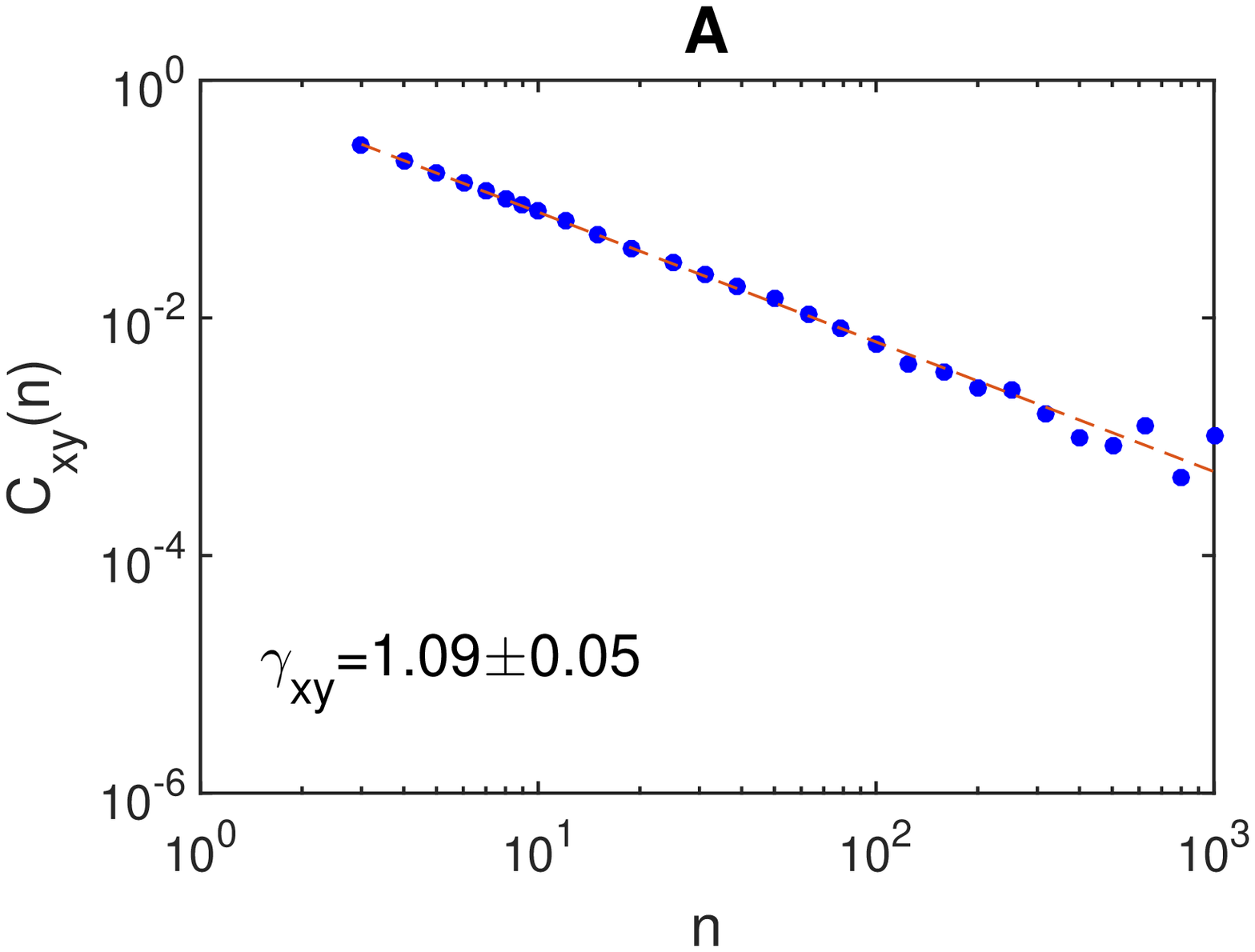}
    }
    \hspace*{-0.9em}
    \subfigure{
        \centering
        \includegraphics[width=0.32\textwidth]{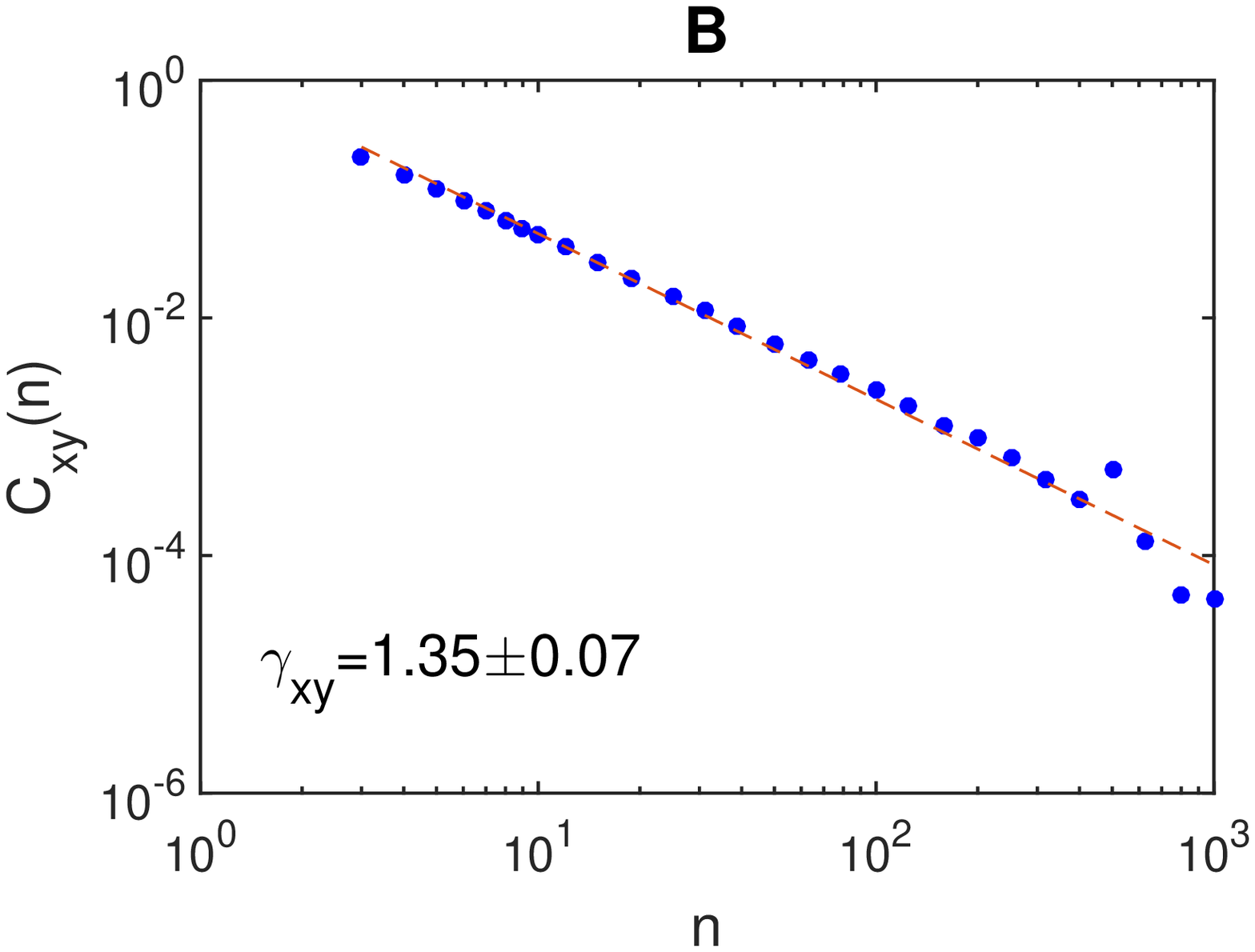}
    }
    \hspace*{-0.9em}
    \subfigure{
        \centering
        \includegraphics[width=0.32\textwidth]{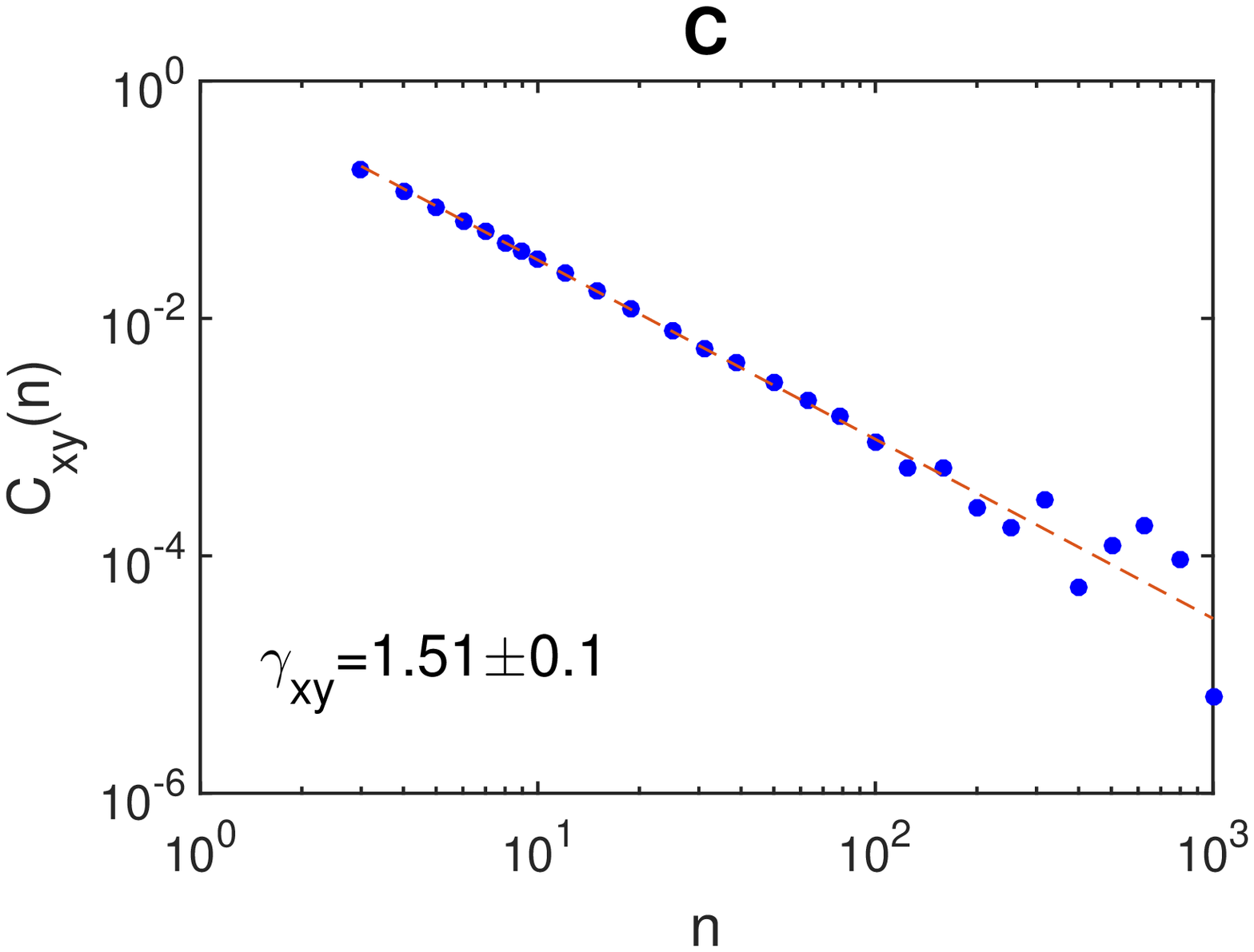}
    }
\caption{The top row represents our numerical results obtained for two coupled self-affine surfaces $h_x$ and $h_y$ of size $256\times 256$. The power-law exponents of the underlying fields $\{x_{\mathbf{i}}\}$ \& $\{y_{\mathbf{i}}\}$ are $\gamma_{xx}=0.7$, $\gamma_{yy}=1.5$, and $\gamma_{xy}=1$. The middle and the bottom panels are log-log plots of average correlations $C_{xx}$, $C_{yy}$, and $C_{xy}$ for three different cases which is labeled by 'A', 'B', and 'C'. Each case is the average of $100$ samples in a square lattice of $2^{12} \times 2^{12}$. The middle panels show the log-log autocorrelation functions for different desired values of $(\gamma_{xx},\gamma_{yy})$ respectively from left to right equal to $(1.3,1.5)$, $(1.1,1.5)$, and $(1.1,1.3)$. The bottom panels show the cross-correlation function for different desired values of $\gamma_{xy}$ respectively from left to right equal to $1.1$, $1.3$, and $1.5$.  The best fit for each exponent is labeled in each panel.}
\label{fig: two coupled random fields}
\end{figure*}

\section{Two coupled random fields}
\label{sec: Two coupled random fields}

Since, sometimes, we may need to generate two coupled random fields~\cite{JVJPJ}, generalizing this method to random fields would be very useful. Consider two discrete $d$-dimensional Gaussian random fields namely, $\{x_{\mathbf{i}}\}$ and $\{y_{\mathbf{i}}\}$ which are defined over a $d$-dimensional cube of volume $L^d$ with desired correlations $C_{xx}(\mathbf{n}), C_{yy}(\mathbf{n})$, and $C_{xy}(\mathbf{n})$ where $\mathbf{i}\equiv(i_1,i_2,\dots,i_d)$ and $\mathbf{n}\equiv(n_1,n_2,\dots,n_d)$. By the definition, the autospectral densities are
\begin{subequations}
\begin{align}
\label{eq: 2d SDxx}
S_{xx}(\mathbf{q}) &=\langle x_{\mathbf{q}} x_{-\mathbf{q}} \rangle \\
\label{eq: 2d SDyy}
S_{yy}(\mathbf{q}) &=\langle y_{\mathbf{q}} y_{-\mathbf{q}} \rangle,
\end{align}
and cross-spectral density is
\begin{align}
\label{eq: 2d SDxy}
S_{xy}(\mathbf{q}) &=\langle x_{\mathbf{q}} y_{-\mathbf{q}} \rangle,
\end{align}
\end{subequations}
where $\{x_{\mathbf{q}}\}$ and $\{y_{\mathbf{q}}\}$ are the Fourier transforms of $\{x_{\mathbf{i}}\}$ and $\{y_{\mathbf{i}}\}$ respectively \cite{EV}. Assuming that the fields are homogeneous and isotropic, the correlations $C_{xx}(\mathbf{n}), C_{yy}(\mathbf{n})$ \& $C_{xy}(\mathbf{n})$ only depend on $n=|\mathbf{n}|$ and the spectral densities $S_{xx}(\mathbf{q}), S_{yy}(\mathbf{q})$ \& $S_{xy}(\mathbf{q})$ which are the Fourier transforms of corresponding correlation functions only depend on $q=|\mathbf{q}|$. This property, the dependency of spectral densities to only one parameter $(q)$, would allow us to apply the procedure of stochastic processes to random fields. In other words, for generating two coupled random fields which are homogeneous and isotropic, the path taken is similar to that of two coupled sequences. Hence, from this point on, we focus on homogeneous and isotropic Gaussian random fields.

The goal here is to construct two coupled random fields $\{x_{\mathbf{i}}\}$ and $\{y_{\mathbf{i}}\}$ with the desired correlations starting from two independent and uncorrelated Gaussian random fields $\{u_{\mathbf{i}}\}$ and $\{v_{\mathbf{i}}\}$. To do so, we combine $\{u_{\mathbf{i}}\}$ and $\{v_{\mathbf{i}}\}$ in the Fourier space as
\begin{equation}
\label{eq: xq and yq in 2d}
\begin{aligned}
x_{\mathbf{q}} &=A_{\mathbf{q}}\, u_{\mathbf{q}}+B_{\mathbf{q}}\, v_{\mathbf{q}} \\
y_{\mathbf{q}} &=C_{\mathbf{q}}\, u_{\mathbf{q}}+D_{\mathbf{q}}\, v_{\mathbf{q}},
\end{aligned}
\end{equation}\\
where $\{u_{\mathbf{q}}\}$ and $\{v_{\mathbf{q}}\}$ are the Fourier transforms of $\{u_{\mathbf{i}}\}$ and $\{v_{\mathbf{i}}\}$ respectively. Since $x_{\mathbf{q}}$ and $y_{\mathbf{q}}$ should satisfy Eqs.~(\ref{eq: 2d SDxx})-(\ref{eq: 2d SDxy}), finally, the functional forms of the above coefficients in terms of spectral densities would be obtained.

Because of dealing with homogeneous and isotropic random fields, we assume that the coefficients $A_{\mathbf{q}}, B_{\mathbf{q}}, C_{\mathbf{q}}$, and $D_{\mathbf{q}}$ are only functions of $q=|{\mathbf{q}}|$. Moreover, we assume that each coefficients has the property that its value at $(-q)$ is equal to its conjugate at $(q)$, similar to Eq.~(\ref{eq: ABCD}). These two assumptions ultimately leads to equations
\begin{subequations}
 \begin{align}
   S_{xx}(q) &= |A_q|^2 + |B_q|^2 \\
   S_{yy}(q) &= |C_q|^2 + |D_q|^2 \\
   S_{xy}(q) &= A_q C_q^* + B_q D_q^*.
  \end{align}
\end{subequations}
The above equations are exactly similar to those in Eqs.~(\ref{eq: SDxx and Coefficients 2}-\ref{eq: SDxy and Coefficients 2}) and so the coefficients $A_{\mathbf{q}}, B_{\mathbf{q}}, C_{\mathbf{q}}$, and $D_{\mathbf{q}}$ are given by the Eqs.~(\ref{eq: Coefficients 1}) and (\ref{eq: Coefficients 2}). As a result, one can generate two Gaussian random fields which are homogeneous and isotropic with the algorithm proposed for two Gaussian random sequences in section \ref{sec: algorithm}.

The question that may arise here is; why are the coefficients and the algorithm similar for random sequences and random fields? The answer to this question lies in the assumption of homogeneity and isotropicity for random fields. To be more precise, when the random fields are assumed to be homogeneous and isotropic, the correlation functions and spectral densities become functions of only one variable, $C(\mathbf{n})=C(n)$ and $S(\mathbf{q})=S(q)$. Our algorithm for generating two coupled random fields is similar to when generating two coupled random sequences. However, the only difference is related to the form of spectral densities. For instance, consider a special case of $d$-dimensional random fields with power-law correlations. Makse et.al. showed that the spectral density of a $d$-dimensional random field with the power-law exponent $\gamma$ has the following form~\cite{MHSS96}
\begin{equation}
\label{eq: d-dimensional power law spectral density}
S(\mathbf{q})=\frac{2\pi^{d/2}}{\Gamma(\beta_d+1)}\left(\frac{q}{2}\right)^{\beta_d} K_{\beta_d}(q),
\end{equation}
where $K_{\beta_d}(q)$ is the modified Bessel function of order $\beta_d=(\gamma-d)/2$, $\Gamma$ is the gamma function. Note that, $q_i=2\pi m_i/L$, $-L/2\le m_i \le L/2$, $i=1, \dots, d$. Equation.~(\ref{eq: d-dimensional power law spectral density}) shows explicitly the dependency on the dimension of space, $d$.

In order to test our algorithm, two-dimensional random fields with power-law correlations is considered. Before proceeding to test the algorithm it is worthy to say that, given a power-law correlated random field $\{x_{i,j}\}$ with $0<\gamma<2$, one can construct a two-dimensional FBM by defining~\cite{MHSS96}
\begin{equation}
h_x(s,t) = \sum_{i=1}^{s} x_{i,t} + \sum_{j=1}^{t} x_{s,j},
\end{equation}
where in the physics literature it is also called a self-affine surface. The top row of Fig.~\ref{fig: two coupled random fields} is our numerical representation of two self-affine surfaces $h_x$ and $h_y$  of size $256\times 256$ with power-law coupling. The power-law exponents of the underlying fields are $\gamma_{xx}=0.7$, $\gamma_{yy}=1.5$, and $\gamma_{xy}=1$. The middle and the bottom panels of Fig.~\ref{fig: two coupled random fields} show the results for autocorrelation functions together with the cross-correlation function for the square lattice of size $2^{12} \times 2^{12}$. Note that these plots are done by considering three various cases 'A', 'B', and 'C'. \\

\section{Conclusion}

We developed an algorithm based on the modified version of the Fourier filtering method to generate two coupled stochastic processes. This algorithm is general because of leaving us free in choosing any desired correlations either \emph{autocorrelations} or \emph{cross-correlation}. Moreover, the method proposed in this work can also be used to generate two coupled random fields of any dimensionality with any sort of correlations. The method has been applied to one and two-dimensional models namely times series and rough surfaces. A random media is a three dimensional application of the present study.

\section{Acknowledgement}
The authors would like to thank Dr. Soheil Vasheghani Farahani, Dr. Mohammad Reza Dehghani, and Dr. Hossein Bayani who assisted in editing the manuscript.

\bibliography{GCCTS}

\begin{thebibliography}{33}%
\makeatletter
\providecommand \@ifxundefined [1]{%
 \@ifx{#1\undefined}
}%
\providecommand \@ifnum [1]{%
 \ifnum #1\expandafter \@firstoftwo
 \else \expandafter \@secondoftwo
 \fi
}%
\providecommand \@ifx [1]{%
 \ifx #1\expandafter \@firstoftwo
 \else \expandafter \@secondoftwo
 \fi
}%
\providecommand \natexlab [1]{#1}%
\providecommand \enquote  [1]{``#1''}%
\providecommand \bibnamefont  [1]{#1}%
\providecommand \bibfnamefont [1]{#1}%
\providecommand \citenamefont [1]{#1}%
\providecommand \href@noop [0]{\@secondoftwo}%
\providecommand \href [0]{\begingroup \@sanitize@url \@href}%
\providecommand \@href[1]{\@@startlink{#1}\@@href}%
\providecommand \@@href[1]{\endgroup#1\@@endlink}%
\providecommand \@sanitize@url [0]{\catcode `\\12\catcode `\$12\catcode
  `\&12\catcode `\#12\catcode `\^12\catcode `\_12\catcode `\%12\relax}%
\providecommand \@@startlink[1]{}%
\providecommand \@@endlink[0]{}%
\providecommand \url  [0]{\begingroup\@sanitize@url \@url }%
\providecommand \@url [1]{\endgroup\@href {#1}{\urlprefix }}%
\providecommand \urlprefix  [0]{URL }%
\providecommand \Eprint [0]{\href }%
\providecommand \doibase [0]{http://dx.doi.org/}%
\providecommand \selectlanguage [0]{\@gobble}%
\providecommand \bibinfo  [0]{\@secondoftwo}%
\providecommand \bibfield  [0]{\@secondoftwo}%
\providecommand \translation [1]{[#1]}%
\providecommand \BibitemOpen [0]{}%
\providecommand \bibitemStop [0]{}%
\providecommand \bibitemNoStop [0]{.\EOS\space}%
\providecommand \EOS [0]{\spacefactor3000\relax}%
\providecommand \BibitemShut  [1]{\csname bibitem#1\endcsname}%
\let\auto@bib@innerbib\@empty
\bibitem [{\citenamefont {Horn}\ and\ \citenamefont {Johnson}(2013)}]{HJ2013}%
  \BibitemOpen
  \bibfield  {author} {\bibinfo {author} {\bibfnamefont {R.~A.}\ \bibnamefont
  {Horn}}\ and\ \bibinfo {author} {\bibfnamefont {C.~R.}\ \bibnamefont
  {Johnson}},\ }\href@noop {} {\emph {\bibinfo {title} {Matrix analysis}}},\
  \bibinfo {edition} {2nd}\ ed.\ (\bibinfo  {publisher} {Cambridge university
  press},\ \bibinfo {year} {2013})\BibitemShut {NoStop}%
\bibitem [{\citenamefont {Golub}\ and\ \citenamefont
  {Van~Loan}(1996)}]{GL1996}%
  \BibitemOpen
  \bibfield  {author} {\bibinfo {author} {\bibfnamefont {G.~H.}\ \bibnamefont
  {Golub}}\ and\ \bibinfo {author} {\bibfnamefont {C.~F.}\ \bibnamefont
  {Van~Loan}},\ }\href@noop {} {\emph {\bibinfo {title} {Matrix
  computations}}},\ \bibinfo {edition} {3rd}\ ed.\ (\bibinfo  {publisher}
  {Johns Hopkins University Press},\ \bibinfo {year} {1996})\BibitemShut
  {NoStop}%
\bibitem [{\citenamefont {Peitgen}\ and\ \citenamefont {Saupe}(1988)}]{SFI}%
  \BibitemOpen
  \bibinfo {editor} {\bibfnamefont {H.-O.}\ \bibnamefont {Peitgen}}\ and\
  \bibinfo {editor} {\bibfnamefont {D.}~\bibnamefont {Saupe}},\ eds.,\ \enquote
  {\bibinfo {title} {The science of fractal images},}\ \ (\bibinfo  {publisher}
  {Springer},\ \bibinfo {year} {1988})\BibitemShut {NoStop}%
\bibitem [{\citenamefont {{Peng}}\ \emph {et~al.}(1991)\citenamefont {{Peng}},
  \citenamefont {{Havlin}}, \citenamefont {{Schwartz}},\ and\ \citenamefont
  {{Stanley}}}]{PHSS91}%
  \BibitemOpen
  \bibfield  {author} {\bibinfo {author} {\bibfnamefont {C.-K.}\ \bibnamefont
  {{Peng}}}, \bibinfo {author} {\bibfnamefont {S.}~\bibnamefont {{Havlin}}},
  \bibinfo {author} {\bibfnamefont {M.}~\bibnamefont {{Schwartz}}}, \ and\
  \bibinfo {author} {\bibfnamefont {H.~E.}\ \bibnamefont {{Stanley}}},\ }\href
  {\doibase 10.1103/PhysRevA.44.R2239} {\bibfield  {journal} {\bibinfo
  {journal} {Phys. Rev. A}\ }\textbf {\bibinfo {volume} {44}},\ \bibinfo
  {pages} {2239} (\bibinfo {year} {1991})}\BibitemShut {NoStop}%
\bibitem [{\citenamefont {{Prakash}}\ \emph {et~al.}(1992)\citenamefont
  {{Prakash}}, \citenamefont {{Havlin}}, \citenamefont {{Schwartz}},\ and\
  \citenamefont {{Stanley}}}]{PHSS92}%
  \BibitemOpen
  \bibfield  {author} {\bibinfo {author} {\bibfnamefont {S.}~\bibnamefont
  {{Prakash}}}, \bibinfo {author} {\bibfnamefont {S.}~\bibnamefont {{Havlin}}},
  \bibinfo {author} {\bibfnamefont {M.}~\bibnamefont {{Schwartz}}}, \ and\
  \bibinfo {author} {\bibfnamefont {H.~E.}\ \bibnamefont {{Stanley}}},\ }\href
  {\doibase 10.1103/PhysRevA.46.R1724} {\bibfield  {journal} {\bibinfo
  {journal} {Phys. Rev. A}\ }\textbf {\bibinfo {volume} {46}},\ \bibinfo
  {pages} {1724} (\bibinfo {year} {1992})}\BibitemShut {NoStop}%
\bibitem [{\citenamefont {{Makse}}\ \emph {et~al.}(1996)\citenamefont
  {{Makse}}, \citenamefont {{Havlin}}, \citenamefont {{Schwartz}},\ and\
  \citenamefont {{Stanley}}}]{MHSS96}%
  \BibitemOpen
  \bibfield  {author} {\bibinfo {author} {\bibfnamefont {H.~A.}\ \bibnamefont
  {{Makse}}}, \bibinfo {author} {\bibfnamefont {S.}~\bibnamefont {{Havlin}}},
  \bibinfo {author} {\bibfnamefont {M.}~\bibnamefont {{Schwartz}}}, \ and\
  \bibinfo {author} {\bibfnamefont {H.~E.}\ \bibnamefont {{Stanley}}},\
  }\href@noop {} {\bibfield  {journal} {\bibinfo  {journal} {Phys. Rev. E}\
  }\textbf {\bibinfo {volume} {53}},\ \bibinfo {pages} {5445} (\bibinfo {year}
  {1996})}\BibitemShut {NoStop}%
\bibitem [{\citenamefont {Wood}\ and\ \citenamefont {Chan}(1994)}]{WC1994}%
  \BibitemOpen
  \bibfield  {author} {\bibinfo {author} {\bibfnamefont {A.~T.~A.}\
  \bibnamefont {Wood}}\ and\ \bibinfo {author} {\bibfnamefont {G.}~\bibnamefont
  {Chan}},\ }\href@noop {} {\bibfield  {journal} {\bibinfo  {journal} {Journal
  of Computational and Graphical Statistics}\ }\textbf {\bibinfo {volume}
  {3}},\ \bibinfo {pages} {409} (\bibinfo {year} {1994})}\BibitemShut {NoStop}%
\bibitem [{\citenamefont {Dietrich}\ and\ \citenamefont {Newsam}(1997)}]{DN97}%
  \BibitemOpen
  \bibfield  {author} {\bibinfo {author} {\bibfnamefont {C.~R.}\ \bibnamefont
  {Dietrich}}\ and\ \bibinfo {author} {\bibfnamefont {G.~N.}\ \bibnamefont
  {Newsam}},\ }\href@noop {} {\bibfield  {journal} {\bibinfo  {journal} {SIAM
  Journal on Scientific Computing}\ }\textbf {\bibinfo {volume} {18}},\
  \bibinfo {pages} {1088} (\bibinfo {year} {1997})}\BibitemShut {NoStop}%
\bibitem [{\citenamefont {Lord}\ \emph {et~al.}(2014)\citenamefont {Lord},
  \citenamefont {Powell},\ and\ \citenamefont {Shardlow}}]{LPS2014}%
  \BibitemOpen
  \bibfield  {author} {\bibinfo {author} {\bibfnamefont {G.~J.}\ \bibnamefont
  {Lord}}, \bibinfo {author} {\bibfnamefont {C.~E.}\ \bibnamefont {Powell}}, \
  and\ \bibinfo {author} {\bibfnamefont {T.}~\bibnamefont {Shardlow}},\
  }\href@noop {} {\emph {\bibinfo {title} {An Introduction to Computational
  Stochastic PDEs}}}\ (\bibinfo  {publisher} {Cambridge University Press},\
  \bibinfo {year} {2014})\BibitemShut {NoStop}%
\bibitem [{\citenamefont {Podobnik}\ and\ \citenamefont {Stanley}(2008)}]{ps}%
  \BibitemOpen
  \bibfield  {author} {\bibinfo {author} {\bibfnamefont {B.}~\bibnamefont
  {Podobnik}}\ and\ \bibinfo {author} {\bibfnamefont {H.~E.}\ \bibnamefont
  {Stanley}},\ }\href {\doibase 10.1103/PhysRevLett.100.084102} {\bibfield
  {journal} {\bibinfo  {journal} {Phys. Rev. Lett.}\ }\textbf {\bibinfo
  {volume} {100}},\ \bibinfo {pages} {084102} (\bibinfo {year}
  {2008})}\BibitemShut {NoStop}%
\bibitem [{\citenamefont {Shadkhoo}\ and\ \citenamefont
  {Jafari}(2009)}]{Shadkhoo}%
  \BibitemOpen
  \bibfield  {author} {\bibinfo {author} {\bibfnamefont {S.}~\bibnamefont
  {Shadkhoo}}\ and\ \bibinfo {author} {\bibfnamefont {G.~R.}\ \bibnamefont
  {Jafari}},\ }\href {\doibase 10.1140/epjb/e2009-00402-2} {\bibfield
  {journal} {\bibinfo  {journal} {EPJ B}\ }\textbf {\bibinfo {volume} {72}},\
  \bibinfo {pages} {679} (\bibinfo {year} {2009})}\BibitemShut {NoStop}%
\bibitem [{\citenamefont {O\ifmmode \acute{s}\else
  \'{s}\fi{}wie\ifmmode~\mbox{\c{}}\else \c{}\fi{}cimka}\ \emph
  {et~al.}(2014)\citenamefont {O\ifmmode \acute{s}\else
  \'{s}\fi{}wie\ifmmode~\mbox{\c{}}\else \c{}\fi{}cimka}, \citenamefont
  {Dro\ifmmode \dot{z}\else \.{z}\fi{}d\ifmmode~\dot{z}\else \.{z}\fi{}},
  \citenamefont {Forczek}, \citenamefont {Jadach},\ and\ \citenamefont
  {Kwapie\ifmmode~\acute{n}\else \'{n}\fi{}}}]{DXA2014}%
  \BibitemOpen
  \bibfield  {author} {\bibinfo {author} {\bibfnamefont {P.}~\bibnamefont
  {O\ifmmode \acute{s}\else \'{s}\fi{}wie\ifmmode~\mbox{\c{}}\else
  \c{}\fi{}cimka}}, \bibinfo {author} {\bibfnamefont {S.}~\bibnamefont
  {Dro\ifmmode \dot{z}\else \.{z}\fi{}d\ifmmode~\dot{z}\else \.{z}\fi{}}},
  \bibinfo {author} {\bibfnamefont {M.}~\bibnamefont {Forczek}}, \bibinfo
  {author} {\bibfnamefont {S.}~\bibnamefont {Jadach}}, \ and\ \bibinfo {author}
  {\bibfnamefont {J.}~\bibnamefont {Kwapie\ifmmode~\acute{n}\else
  \'{n}\fi{}}},\ }\href {\doibase 10.1103/PhysRevE.89.023305} {\bibfield
  {journal} {\bibinfo  {journal} {Phys. Rev. E}\ }\textbf {\bibinfo {volume}
  {89}},\ \bibinfo {pages} {023305} (\bibinfo {year} {2014})}\BibitemShut
  {NoStop}%
\bibitem [{\citenamefont {Qian}\ \emph {et~al.}(2015)\citenamefont {Qian},
  \citenamefont {Liu}, \citenamefont {Jiang}, \citenamefont {Podobnik},
  \citenamefont {Zhou},\ and\ \citenamefont {Stanley}}]{partialDXA}%
  \BibitemOpen
  \bibfield  {author} {\bibinfo {author} {\bibfnamefont {X.-Y.}\ \bibnamefont
  {Qian}}, \bibinfo {author} {\bibfnamefont {Y.-M.}\ \bibnamefont {Liu}},
  \bibinfo {author} {\bibfnamefont {Z.-Q.}\ \bibnamefont {Jiang}}, \bibinfo
  {author} {\bibfnamefont {B.}~\bibnamefont {Podobnik}}, \bibinfo {author}
  {\bibfnamefont {W.-X.}\ \bibnamefont {Zhou}}, \ and\ \bibinfo {author}
  {\bibfnamefont {H.~E.}\ \bibnamefont {Stanley}},\ }\href {\doibase
  10.1103/PhysRevE.91.062816} {\bibfield  {journal} {\bibinfo  {journal} {Phys.
  Rev. E}\ }\textbf {\bibinfo {volume} {91}},\ \bibinfo {pages} {062816}
  (\bibinfo {year} {2015})}\BibitemShut {NoStop}%
\bibitem [{\citenamefont {Hedayatifar}\ \emph {et~al.}(2011)\citenamefont
  {Hedayatifar}, \citenamefont {Vahabi},\ and\ \citenamefont {Jafari}}]{HVJ}%
  \BibitemOpen
  \bibfield  {author} {\bibinfo {author} {\bibfnamefont {L.}~\bibnamefont
  {Hedayatifar}}, \bibinfo {author} {\bibfnamefont {M.}~\bibnamefont {Vahabi}},
  \ and\ \bibinfo {author} {\bibfnamefont {G.~R.}\ \bibnamefont {Jafari}},\
  }\href {\doibase 10.1103/PhysRevE.84.021138} {\bibfield  {journal} {\bibinfo
  {journal} {Phys. Rev. E}\ }\textbf {\bibinfo {volume} {84}},\ \bibinfo
  {pages} {021138} (\bibinfo {year} {2011})}\BibitemShut {NoStop}%
\bibitem [{\citenamefont {Ge}(2008)}]{wavelet}%
  \BibitemOpen
  \bibfield  {author} {\bibinfo {author} {\bibfnamefont {Z.}~\bibnamefont
  {Ge}},\ }\href {\doibase 10.5194/angeo-26-3819-2008} {\bibfield  {journal}
  {\bibinfo  {journal} {Annales Geophysicae}\ }\textbf {\bibinfo {volume}
  {26}},\ \bibinfo {pages} {3819} (\bibinfo {year} {2008})}\BibitemShut
  {NoStop}%
\bibitem [{\citenamefont {Dashtian}\ \emph {et~al.}(2011)\citenamefont
  {Dashtian}, \citenamefont {Jafari}, \citenamefont {Koohi~Lai}, \citenamefont
  {Masihi},\ and\ \citenamefont {Sahimi}}]{djkms}%
  \BibitemOpen
  \bibfield  {author} {\bibinfo {author} {\bibfnamefont {H.}~\bibnamefont
  {Dashtian}}, \bibinfo {author} {\bibfnamefont {G.}~\bibnamefont {Jafari}},
  \bibinfo {author} {\bibfnamefont {Z.}~\bibnamefont {Koohi~Lai}}, \bibinfo
  {author} {\bibfnamefont {M.}~\bibnamefont {Masihi}}, \ and\ \bibinfo {author}
  {\bibfnamefont {M.}~\bibnamefont {Sahimi}},\ }\href {\doibase
  10.1007/s11242-011-9794-x} {\bibfield  {journal} {\bibinfo  {journal}
  {Transport in Porous Media}\ }\textbf {\bibinfo {volume} {90}},\ \bibinfo
  {pages} {445} (\bibinfo {year} {2011})}\BibitemShut {NoStop}%
\bibitem [{\citenamefont {Laloux}\ \emph {et~al.}(1999)\citenamefont {Laloux},
  \citenamefont {Cizeau}, \citenamefont {Bouchaud},\ and\ \citenamefont
  {Potters}}]{lcbp}%
  \BibitemOpen
  \bibfield  {author} {\bibinfo {author} {\bibfnamefont {L.}~\bibnamefont
  {Laloux}}, \bibinfo {author} {\bibfnamefont {P.}~\bibnamefont {Cizeau}},
  \bibinfo {author} {\bibfnamefont {J.-P.}\ \bibnamefont {Bouchaud}}, \ and\
  \bibinfo {author} {\bibfnamefont {M.}~\bibnamefont {Potters}},\ }\href
  {\doibase 10.1103/PhysRevLett.83.1467} {\bibfield  {journal} {\bibinfo
  {journal} {Phys. Rev. Lett.}\ }\textbf {\bibinfo {volume} {83}},\ \bibinfo
  {pages} {1467} (\bibinfo {year} {1999})}\BibitemShut {NoStop}%
\bibitem [{\citenamefont {Plerou}\ \emph {et~al.}(2000)\citenamefont {Plerou},
  \citenamefont {Gopikrishnan}, \citenamefont {Rosenow}, \citenamefont
  {Amaral},\ and\ \citenamefont {Stanley}}]{pgras2}%
  \BibitemOpen
  \bibfield  {author} {\bibinfo {author} {\bibfnamefont {V.}~\bibnamefont
  {Plerou}}, \bibinfo {author} {\bibfnamefont {P.}~\bibnamefont
  {Gopikrishnan}}, \bibinfo {author} {\bibfnamefont {B.}~\bibnamefont
  {Rosenow}}, \bibinfo {author} {\bibfnamefont {L.}~\bibnamefont {Amaral}}, \
  and\ \bibinfo {author} {\bibfnamefont {H.}~\bibnamefont {Stanley}},\ }\href
  {\doibase 10.1016/S0378-4371(00)00376-9} {\bibfield  {journal} {\bibinfo
  {journal} {Phys. A}\ }\textbf {\bibinfo {volume} {287}},\ \bibinfo {pages}
  {374 } (\bibinfo {year} {2000})}\BibitemShut {NoStop}%
\bibitem [{\citenamefont {Plerou}\ \emph {et~al.}(2002)\citenamefont {Plerou},
  \citenamefont {Gopikrishnan}, \citenamefont {Rosenow}, \citenamefont
  {Amaral}, \citenamefont {Guhr},\ and\ \citenamefont {Stanley}}]{pgrags}%
  \BibitemOpen
  \bibfield  {author} {\bibinfo {author} {\bibfnamefont {V.}~\bibnamefont
  {Plerou}}, \bibinfo {author} {\bibfnamefont {P.}~\bibnamefont
  {Gopikrishnan}}, \bibinfo {author} {\bibfnamefont {B.}~\bibnamefont
  {Rosenow}}, \bibinfo {author} {\bibfnamefont {L.~A.~N.}\ \bibnamefont
  {Amaral}}, \bibinfo {author} {\bibfnamefont {T.}~\bibnamefont {Guhr}}, \ and\
  \bibinfo {author} {\bibfnamefont {H.~E.}\ \bibnamefont {Stanley}},\ }\href
  {\doibase 10.1103/PhysRevE.65.066126} {\bibfield  {journal} {\bibinfo
  {journal} {Phys. Rev. E}\ }\textbf {\bibinfo {volume} {65}},\ \bibinfo
  {pages} {066126} (\bibinfo {year} {2002})}\BibitemShut {NoStop}%
\bibitem [{\citenamefont {Jamali}\ \emph
  {et~al.}(2015{\natexlab{a}})\citenamefont {Jamali}, \citenamefont {Jafari},\
  and\ \citenamefont {Farahani}}]{JJF}%
  \BibitemOpen
  \bibfield  {author} {\bibinfo {author} {\bibfnamefont {T.}~\bibnamefont
  {Jamali}}, \bibinfo {author} {\bibfnamefont {G.}~\bibnamefont {Jafari}}, \
  and\ \bibinfo {author} {\bibfnamefont {S.~V.}\ \bibnamefont {Farahani}},\
  }\href@noop {} {\  (\bibinfo {year} {2015}{\natexlab{a}})},\ \Eprint
  {http://arxiv.org/abs/1505.03336} {arXiv:1505.03336 [physics.data-an]}
  \BibitemShut {NoStop}%
\bibitem [{\citenamefont {Mehraban}\ \emph {et~al.}(2013)\citenamefont
  {Mehraban}, \citenamefont {Shirazi}, \citenamefont {Zamani},\ and\
  \citenamefont {Jafari}}]{Mehraban}%
  \BibitemOpen
  \bibfield  {author} {\bibinfo {author} {\bibfnamefont {S.}~\bibnamefont
  {Mehraban}}, \bibinfo {author} {\bibfnamefont {A.~H.}\ \bibnamefont
  {Shirazi}}, \bibinfo {author} {\bibfnamefont {M.}~\bibnamefont {Zamani}}, \
  and\ \bibinfo {author} {\bibfnamefont {G.~R.}\ \bibnamefont {Jafari}},\
  }\href {http://stacks.iop.org/0295-5075/103/i=5/a=50011} {\bibfield
  {journal} {\bibinfo  {journal} {EPL (Europhysics Letters)}\ }\textbf
  {\bibinfo {volume} {103}},\ \bibinfo {pages} {50011} (\bibinfo {year}
  {2013})}\BibitemShut {NoStop}%
\bibitem [{\citenamefont {{Podobnik}}\ \emph {et~al.}(2007)\citenamefont
  {{Podobnik}}, \citenamefont {{Fu}}, \citenamefont {{Stanley}},\ and\
  \citenamefont {{Ivanov}}}]{PFSI}%
  \BibitemOpen
  \bibfield  {author} {\bibinfo {author} {\bibfnamefont {B.}~\bibnamefont
  {{Podobnik}}}, \bibinfo {author} {\bibfnamefont {D.~F.}\ \bibnamefont
  {{Fu}}}, \bibinfo {author} {\bibfnamefont {H.~E.}\ \bibnamefont {{Stanley}}},
  \ and\ \bibinfo {author} {\bibfnamefont {P.~C.}\ \bibnamefont {{Ivanov}}},\
  }\href {\doibase 10.1140/epjb/e2007-00089-3} {\bibfield  {journal} {\bibinfo
  {journal} {European Physical Journal B}\ }\textbf {\bibinfo {volume} {56}},\
  \bibinfo {pages} {47} (\bibinfo {year} {2007})}\BibitemShut {NoStop}%
\bibitem [{\citenamefont {{Podobnik}}\ \emph {et~al.}(2008)\citenamefont
  {{Podobnik}}, \citenamefont {{Horvatic}}, \citenamefont {{Lam Ng}},
  \citenamefont {{Eugene Stanley}},\ and\ \citenamefont {{Ivanov}}}]{PHLEI}%
  \BibitemOpen
  \bibfield  {author} {\bibinfo {author} {\bibfnamefont {B.}~\bibnamefont
  {{Podobnik}}}, \bibinfo {author} {\bibfnamefont {D.}~\bibnamefont
  {{Horvatic}}}, \bibinfo {author} {\bibfnamefont {A.}~\bibnamefont {{Lam
  Ng}}}, \bibinfo {author} {\bibfnamefont {H.}~\bibnamefont {{Eugene
  Stanley}}}, \ and\ \bibinfo {author} {\bibfnamefont {P.~C.}\ \bibnamefont
  {{Ivanov}}},\ }\href@noop {} {\bibfield  {journal} {\bibinfo  {journal}
  {Physica A Statistical Mechanics and its Applications}\ }\textbf {\bibinfo
  {volume} {387}},\ \bibinfo {pages} {3954} (\bibinfo {year}
  {2008})}\BibitemShut {NoStop}%
\bibitem [{\citenamefont {Bendat}\ and\ \citenamefont
  {Piersol}(2000)}]{Bendat}%
  \BibitemOpen
  \bibfield  {author} {\bibinfo {author} {\bibfnamefont {J.~S.}\ \bibnamefont
  {Bendat}}\ and\ \bibinfo {author} {\bibfnamefont {A.~G.}\ \bibnamefont
  {Piersol}},\ }\href@noop {} {\emph {\bibinfo {title} {Random Data: Analysis
  and Measurement Procedures}}},\ \bibinfo {edition} {3rd}\ ed.\ (\bibinfo
  {publisher} {John Wiley \& Sons, Inc.},\ \bibinfo {address} {New York, NY,
  USA},\ \bibinfo {year} {2000})\BibitemShut {NoStop}%
\bibitem [{\citenamefont {Wiener}(1964)}]{Wiener}%
  \BibitemOpen
  \bibfield  {author} {\bibinfo {author} {\bibfnamefont {N.}~\bibnamefont
  {Wiener}},\ }\href@noop {} {\emph {\bibinfo {title} {Time Series}}}\
  (\bibinfo  {publisher} {M.I.T. Press},\ \bibinfo {year} {1964})\BibitemShut
  {NoStop}%
\bibitem [{\citenamefont {Larson}\ and\ \citenamefont {Falvo}(2008)}]{LF}%
  \BibitemOpen
  \bibfield  {author} {\bibinfo {author} {\bibfnamefont {R.}~\bibnamefont
  {Larson}}\ and\ \bibinfo {author} {\bibfnamefont {D.~C.}\ \bibnamefont
  {Falvo}},\ }\href@noop {} {\emph {\bibinfo {title} {Elementary Linear
  Algebra}}},\ \bibinfo {edition} {6th}\ ed.\ (\bibinfo  {publisher} {Brooks
  Cole},\ \bibinfo {year} {2008})\BibitemShut {NoStop}%
\bibitem [{\citenamefont {Scharnhorst}(2001)}]{Scharnhorst}%
  \BibitemOpen
  \bibfield  {author} {\bibinfo {author} {\bibfnamefont {K.}~\bibnamefont
  {Scharnhorst}},\ }\href {\doibase 10.1023/A:1012692601098} {\bibfield
  {journal} {\bibinfo  {journal} {Acta Applicandae Mathematica}\ }\textbf
  {\bibinfo {volume} {69}},\ \bibinfo {pages} {95} (\bibinfo {year}
  {2001})}\BibitemShut {NoStop}%
\bibitem [{\citenamefont {{Makse}}\ \emph {et~al.}(1995)\citenamefont
  {{Makse}}, \citenamefont {{Havlin}}, \citenamefont {{Stanley}},\ and\
  \citenamefont {{Schwartz}}}]{MHSS95}%
  \BibitemOpen
  \bibfield  {author} {\bibinfo {author} {\bibfnamefont {H.}~\bibnamefont
  {{Makse}}}, \bibinfo {author} {\bibfnamefont {S.}~\bibnamefont {{Havlin}}},
  \bibinfo {author} {\bibfnamefont {H.~E.}\ \bibnamefont {{Stanley}}}, \ and\
  \bibinfo {author} {\bibfnamefont {M.}~\bibnamefont {{Schwartz}}},\ }\href
  {\doibase 10.1016/0960-0779(95)80035-F} {\bibfield  {journal} {\bibinfo
  {journal} {Chaos Solitons and Fractals}\ }\textbf {\bibinfo {volume} {6}},\
  \bibinfo {pages} {295} (\bibinfo {year} {1995})}\BibitemShut {NoStop}%
\bibitem [{\citenamefont {{Hamzehpour}}\ and\ \citenamefont
  {{Sahimi}}(2006)}]{HS}%
  \BibitemOpen
  \bibfield  {author} {\bibinfo {author} {\bibfnamefont {H.}~\bibnamefont
  {{Hamzehpour}}}\ and\ \bibinfo {author} {\bibfnamefont {M.}~\bibnamefont
  {{Sahimi}}},\ }\href {\doibase 10.1103/PhysRevE.73.056121} {\bibfield
  {journal} {\bibinfo  {journal} {Phys. Rev. E}\ }\textbf {\bibinfo {volume}
  {73}},\ \bibinfo {eid} {056121} (\bibinfo {year} {2006})}\BibitemShut
  {NoStop}%
\bibitem [{\citenamefont {Vahabi}\ and\ \citenamefont {Jafari}(2012)}]{VJ}%
  \BibitemOpen
  \bibfield  {author} {\bibinfo {author} {\bibfnamefont {M.}~\bibnamefont
  {Vahabi}}\ and\ \bibinfo {author} {\bibfnamefont {G.~R.}\ \bibnamefont
  {Jafari}},\ }\href@noop {} {\bibfield  {journal} {\bibinfo  {journal} {Phys.
  Rev. E}\ }\textbf {\bibinfo {volume} {86}},\ \bibinfo {pages} {066704}
  (\bibinfo {year} {2012})}\BibitemShut {NoStop}%
\bibitem [{\citenamefont {Plerou}\ \emph {et~al.}(2001)\citenamefont {Plerou},
  \citenamefont {Gopikrishnan}, \citenamefont {Gabaix}, \citenamefont
  {Amaral},\ and\ \citenamefont {Stanley}}]{PGGAS}%
  \BibitemOpen
  \bibfield  {author} {\bibinfo {author} {\bibfnamefont {V.}~\bibnamefont
  {Plerou}}, \bibinfo {author} {\bibfnamefont {P.}~\bibnamefont
  {Gopikrishnan}}, \bibinfo {author} {\bibfnamefont {X.}~\bibnamefont
  {Gabaix}}, \bibinfo {author} {\bibfnamefont {L.}~\bibnamefont {Amaral}}, \
  and\ \bibinfo {author} {\bibfnamefont {H.}~\bibnamefont {Stanley}},\ }\href
  {\doibase 10.1088/1469-7688/1/2/308} {\bibfield  {journal} {\bibinfo
  {journal} {Quantitative Finance}\ }\textbf {\bibinfo {volume} {1}},\ \bibinfo
  {pages} {262} (\bibinfo {year} {2001})}\BibitemShut {NoStop}%
\bibitem [{\citenamefont {Jamali}\ \emph
  {et~al.}(2015{\natexlab{b}})\citenamefont {Jamali}, \citenamefont {Farahani},
  \citenamefont {Jannesar}, \citenamefont {Palasantzas},\ and\ \citenamefont
  {Jafari}}]{JVJPJ}%
  \BibitemOpen
  \bibfield  {author} {\bibinfo {author} {\bibfnamefont {T.}~\bibnamefont
  {Jamali}}, \bibinfo {author} {\bibfnamefont {S.~V.}\ \bibnamefont
  {Farahani}}, \bibinfo {author} {\bibfnamefont {M.}~\bibnamefont {Jannesar}},
  \bibinfo {author} {\bibfnamefont {G.}~\bibnamefont {Palasantzas}}, \ and\
  \bibinfo {author} {\bibfnamefont {G.~R.}\ \bibnamefont {Jafari}},\ }\href
  {\doibase http://dx.doi.org/10.1063/1.4919817} {\bibfield  {journal}
  {\bibinfo  {journal} {Journal of Applied Physics}\ }\textbf {\bibinfo
  {volume} {117}},\ \bibinfo {eid} {175308} (\bibinfo {year}
  {2015}{\natexlab{b}})}\BibitemShut {NoStop}%
\bibitem [{\citenamefont {Vanmarcke}(2010)}]{EV}%
  \BibitemOpen
  \bibfield  {author} {\bibinfo {author} {\bibfnamefont {E.}~\bibnamefont
  {Vanmarcke}},\ }\href@noop {} {\emph {\bibinfo {title} {Random Fields:
  Analysis and Synthesis}}}\ (\bibinfo  {publisher} {World Scientific
  Publishing Company; Revised and Expanded New Edition edition},\ \bibinfo
  {year} {2010})\BibitemShut {NoStop}%
\end{thebibliography}%

\end{document}